\documentclass[a4paper,11pt,reqno,english]{smfart}
\usepackage{psfrag}
\usepackage{latexsym}
\usepackage{array}
\usepackage{amssymb}
\theoremstyle{plain}

\usepackage{amsmath,amssymb,amsfonts,graphicx,psfrag,subfig}
\usepackage[latin1]{inputenc}

\usepackage[colorlinks=true, pdfstartview=FitV, linkcolor=blue, citecolor=blue, urlcolor=blue]{hyperref}
\newtheorem{thm}{Theorem}[section]
\newtheorem{prop}[thm]{Proposition}
\newtheorem{lem}[thm]{Lemma}
\newtheorem{df}[thm]{Definition}

\newtheorem{rem}[thm]{Remark}

\def\be#1 {\begin{equation} \label{#1}}
\newcommand{\ee}{\end{equation}}
\def\dem {\noindent {\bf Proof: }}

\newcommand{\mb}{\medskip\noindent}
\newcommand{\gb}{\bigskip\noindent}
\newcommand{\R}{\mathbb R}

\newcommand{\C}{\mathcal C}

\newcommand{\ww}{\langle}
\newcommand{\xx}{\rangle}

\newcommand{\I}{\mathcal I}

\newcommand{\PPP}{\mathrm P}

\def \NN {\mathrm{N}}
\def \TT {\mathrm{T}}

\def \virg {\, , \,\,}
\def \dsp {\displaystyle}
\def \vsp {\vspace{6pt}}

\def\sqw{\hbox{\rlap{\leavevmode\raise.3ex\hbox{$\sqcap$}}$%
\sqcup$}}
\def\findem{\ifmmode\sqw\else{\ifhmode\unskip\fi\nobreak\hfil
\penalty50\hskip1em\null\nobreak\hfil\sqw
\parfillskip=0pt\finalhyphendemerits=0\endgraf}\fi}

\begin{document}

\title[Existence of solutions for second-order differential inclusions]{Existence of solutions for second-order differential inclusions involving proximal normal cones}
\date{\today}

\author{Fr\'ed\'eric Bernicot}
\address{CNRS - Universit\'e Lille 1 \\  Laboratoire de math\'ematiques Paul Painlev\'e \\ 59655 Villeneuve d'Ascq Cedex, France}
\email{frederic.bernicot@math.univ-lille1.fr} \urladdr{http://math.univ-lille1.fr/$\sim$bernicot/}
\author{Juliette Venel}
\address{LAMAV \\ Universit\'e de Valenciennes et du Hainaut-Cambr\'esis Campus \\ Mont Houy 59313 Valenciennes Cedex 9, France}
\email{juliette.venel@univ-valenciennes.fr} \urladdr{http://www.math.u-psud.fr/$\sim$venel/}

\maketitle

\frontmatter 
\begin{abstract} In this work, we prove global existence of solutions for second order differential problems in a general framework. More precisely, we consider second order differential inclusions involving proximal normal cone to a set-valued map. This set-valued map is supposed to take admissible values (so in particular uniformly prox-regular values, which may be non-smooth and non-convex). Moreover we require the solution to satisfy an {\it impact law}, appearing in the description of mechanical systems with inelastic shocks. 
\end{abstract}

\begin{altabstract} Nous prouvons dans ce travail un r\'esultat d'existence globale de solutions pour des inclusions diff\'erentielles du second ordre. Nous consid\'erons plus pr\'ecis\'ement des inclusions diff\'erentielles du second ordre faisant appara\^ itre le c\^ one proximal normal d'un ensemble $C$ d\'ependant du temps. La multifonction $C(\cdot)$ est suppos\'ee admissible (et en particulier prendre des valeurs uniform\'ement prox-r\'eguli\`eres \'eventuellement non r\'eguli\`eres et non convexes). De plus, nous imposons \`a la solution de v\'erifier une {\it loi d'impact}, apparaissant dans la description de syst\`emes m\'ecaniques avec des chocs in\'elastiques. 
\end{altabstract}

\altkeywords{ Inclusions diff\'erentielles du second ordre ; cone proximal normal ; existence de solutions}

\tableofcontents

\section{Introduction}

We consider second order differential inclusions, involving proximal normal cones. These differential inclusions appear in several fields (granular media \cite{Moreau1, Moreau2}, robotics \cite{ER}, \cite{F} and virtual reality \cite{P} ... ).

They describe an evolution problem where the state-variable is submitted to some constraints and so has to live in an {\it admissible set}).

\mb Numerous works deal with the particular case where the admissible set is given by several constraints. Let us first detail them.
We write $C\subset \R^d$ for the closed admissible set defined by constraints $(g_i)_{1\leq i\leq p}$ as follows~:
\be{C:inter} C := \bigcap_{i=1}^p \left\{q\in\R^d,\  g_i(q)\geq 0 \right\}. \ee
The tangent cone to $C$ at $q$ is 
$$ \TT_C(q):= \left\{ u\in \R^d,\ \langle u,\nabla g_i(q)\rangle \geq 0 \quad \forall i\in I(q)\right\}$$
where $I(q)$ is the set of ``active constraints''
$$ I(q):=\{i,\ g_i(q)=0\}.$$
The second order problem is the following one: let $\I$ be a bounded time-interval, $f:\I\times \R^d \rightarrow \R^d$ be a map, find $q\in W^{1,\infty}(\I,\R^d)$ such that $\dot{q}\in BV(\I,\R^d)$ and
\begin{equation} \label{Pb:}
\left\{
\begin{array}{l}
\dsp d \dot{q} + N_C(q) dt \ni f(t, q)dt  \vsp \\
\dsp \forall t\in \I, \quad \dot q(t^+) = \PPP_{\TT_C(q(t))}(\dot q(t^-)) \vsp \\
\dsp q(0)=q_0 \in \textrm{Int}(C) \vsp \\
\dsp \dot q (0)=u_0,
\end{array}
\right.
\end{equation}
where $\textrm{Int}(C)$ is the interior of the set $C$ and $N_C(q)$ is the normal cone defined by
$$ N_C(q):= \left\{ \begin{array}{ll}
                   \{0\} & \textrm{if } q\in \textrm{Int}(C) \vsp \\
                   \left\{-\sum_{i\in I(q)} \lambda_i \nabla g_i(q), \quad \lambda_i\geq 0\right\} & \textrm{if } q\in \partial C \vsp  \\ 
 \emptyset & \textrm{if } q\notin C. \vsp  \end{array} \right. $$
This differential inclusion can be thought as follows: the point $q(t)$, submitted to the external force $f(t,q(t))$, has to live in the set $C$. The differential inclusion $d \dot{q} + N_C(q) dt \ni f(t, q)dt$ does not uniquely define the evolution of the velocity during an impact. To complete the description, we impose the impact law $\dot q(t^+) = \PPP_{\TT_C(q(t))}\dot q(t^-),$
introduced by J.J. Moreau in \cite{Moreau1} and justified by L. Paoli and M. Schatzman in \cite{Paoli-Scha2, Paoli-Scha4} (using a penalty method) for inelastic impacts. This law can be extended with a restitution coefficient to model the elastic shocks.

\mb
The existence of a solution for such second-order problems is still open in a general framework. The first positive results were obtained by M.P.D. Monteiro Marques \cite{Marques}, L. Paoli and M. Schatzman \cite{Paoli-Scha3} in the case of a smooth admissible set (which locally corresponds to the single constraint case $p=1$ in (\ref{C:inter})). The single constraint case is also treated in \cite{MP, DMP} where an additional mass-matrix depending on $q$. The proofs use a time-discretization of (\ref{Pb:}) and rely on the convergence of the approximate solutions. The multi-constraint case with analytical data was then treated by P. Ballard with a different method in \cite{Ballard}. In this paper, a positive result of uniqueness for such problems was also obtained. Then in \cite{Paoli}, an existence result is proved in the case of a non-smooth convex admissible set $C$ (given by multiple constraints). There, the active constraints are supposed to be linearly independent in the following sense: for each $q\in \partial C$, the gradients $(\nabla g_i(q))_{i\in I(q)}$ are supposed to be linearly independent. This assumption is quite strong since it implies that the number of active constraints $|I(q)|$ is always lower than the dimension $d$ (which may fail for some applications). Moreover the impact law is proved under a geometrical assumption:
\be{eq:angle} \forall (i,\ j) \in I(q)^2,\quad  i\neq j,  \qquad \langle \nabla g_i(q), \nabla g_j(q)\rangle \leq 0. \ee
In \cite{Paoli2, Paoli3}, similar results are obtained without requiring the convexity of $C$. Note that excepted in \cite{Ballard}, all these results are local in the following sense: for each initial data $(q_0,u_0)$, there exists a time $\tau=\tau(|u_0|)$ and a solution $q$ of (\ref{Pb:}) on $\I=[0,\tau]$. 

\gb

Recently, time-dependent constraints are considered. More precisely, the constraints $g_i : \I \times \R^d \rightarrow \R$ define the following set-valued map
\be{C:interb} C(\cdot) := \bigcap_{i=1}^p \left\{q\in\R^d,\  g_i(\cdot,q)\geq 0 \right\}. \ee
Then Problem (\ref{Pb:}) becomes: find $q\in W^{1,\infty}(\I,\R^d)$ such that $\dot{q}\in BV(\I,\R^d)$
\begin{equation} \label{Pb:2}
\left\{
\begin{array}{l}
\dsp d \dot{q} + N_{C(t)}(q)dt \ni f(t, q)dt \vsp \\
\dsp \forall t\in \I, \quad \dot q(t^+) = \PPP_{V_{t,q(t)}}\dot q(t^-)\vsp \\
\dsp q(0)=q_0 \in \textrm{Int}(C(0)) \vsp \\
\dsp \dot q (0)=u_0,
\end{array}
\right.
\end{equation}
where $V_{t,q}$ is the set of {\it admissible velocities}: 
\be{def:Cq}  V_{t,q}:=\left\lbrace u,\ \partial_t g_i(t,q) + \langle \nabla_q \, g_i
(t,q),u\rangle\geq 0\ \textrm{ if } \  g_i(t,q)=0 \right\rbrace.\ee
Local existence results were obtained in \cite{Sc3} assuming that $\partial C(\cdot)$ belongs to $C^3(\I \times \R^d)$ (which locally corresponds to a single constraint). Recently in \cite{F-Aline}, a global result was proved by the first author and A. Lefebvre for convex $C^2$ constraints $g_i$ (no regularity on the boundaries $\partial C(t)$ is required). Moreover, Assumption (\ref{eq:angle}) and the independence of the gradients $(\nabla_q g_i(t,q))_{i\in I(t,q)}$ are relaxed to a positive linearly independence: global results are shown under the existence of $\rho,\gamma>0$ such that for all $q\in C(t)$ and nonnegative reals $\lambda_i$
\be{reverse:intro}
\sum_{i\in I_\rho(t,q)} \lambda_i |\nabla_q g_i(t,q)| \leq \gamma \left| \sum_{i\in I_\rho(t,q)} \lambda_i \nabla_q g_i(t,q)\right |, \tag{$R_\rho$}
\ee
where
$$ I_\rho(t,q):=\left\{i,\ g_i(t,q)\leq \rho\right\}.$$

\gb

In this paper, we are interested in local and global existence results for second order problems involving a general set-valued map $C(\cdot)$ (not necessary defined with constraints). This problem was solved in the case of first order differential inclusions called {\it sweeping process}. Let us briefly recall the corresponding context. For a bounded time-interval $\I$, a set-valued map $C:\I \rightrightarrows \R^d$ with nonempty closed values and a map $f: \I \times \R^d \rightrightarrows \R^d$, the associated sweeping process takes the following form: 
\begin{equation} 
\left\{
\begin{array}{l}
 \dsp \frac{dq}{dt}(t) + \NN(C(t),q(t)) \ni f(t,q(t)) \vsp \\
 q(0)=q_0\ ,
\end{array}
\right. \label{sys1}
\end{equation}
with an initial data $q_0\in C(0)$ and where $\NN(C,q)$ denotes the proximal normal cone of $C$ at any point $q$. This kind of evolution problem has been introduced by
 J.J. Moreau in 70's  (see \cite{Moreausweep}) with convex sets $C(t)$. He proposed a time-discretization of (\ref{sys1}) called the \textit{catching-up algorithm}. This scheme was later adapted to the second order problems (\ref{Pb:}) and (\ref{Pb:2}). Later the convexity assumption was weakened by the concept of ``uniform prox-regularity'' (a set $C$ is said to be {\it uniformly prox-regular with constant $\eta$} or {\it $\eta$-prox-regular} if the projection onto $C$ is single-valued and continuous at any point distant at most $\eta$ from $C$). In this framework, the well-posedness of (\ref{sys1}) was proved by G.~Colombo, V.V.~Goncharov in \cite{Colombo}, H.~Benabdellah in \cite{Benab}, L.~Thibault \cite{Thibsweep} and by G. Colombo, M.D.P.~Monteiro Marques in \cite{Monteiro}. The sweeping process problem is still extensively studied but the recent results are not detailed here (see the works of M.~Bounkhel, J.F.~Edmond and L.~Thibault in \cite{Thibnonconv, Thibsweep, Thibrelax, Thibbv} and of the authors in \cite{BV, F-Juliette}).

\mb Note that a major difference exists between the first and the second order differential inclusions. Indeed, even with smooth data the uniqueness does not hold for second order problem such (\ref{Pb:}) or (\ref{Pb:2}), see \cite{Sc} and \cite{Ballard}. The only positive results are proved in \cite{Sc2} for one-dimensional impact problems and in \cite{Ballard} in the context of analytic data. Here we are only interested in existence results for general second order problems.
 
\mb In the case of a general set-valued map $C(\cdot)$, we introduce the following set (which extends the concept of $\TT_C(q)$ and $V_{t,q}$):
\begin{df} \label{def:cone} For a set-valued map $C:\I \rightrightarrows \R^d$, for every $t_0\in \overset{\circ}{\I}$ and $q_0\in C(t)$, the cone of admissible velocities is defined as follows:
$$ \C_{t_0,q_0}:=\left\{v =\lim_{\epsilon \searrow 0} v_\epsilon, \ \textrm{ with } \ v_\epsilon \in \frac{C(t_0+\epsilon)-q_0}{\epsilon}\right\}=\liminf_{\epsilon \searrow 0} \ \frac{C(t_0+\epsilon)-q_0}{\epsilon}.$$ 
\end{df}

\mb 
We assumed that the set-valued map $C(\cdot)$ is Lipschitz continuous on $\I$: there exists $c_0>0$ such that for all $t,s\in \I$
\be{eq:Qlip} d_H(C(t),C(s))\leq c_0|t-s|,\ee
where $d_H$ denotes the Hausdorff distance. \\
We are interested in the following problem: find $q\in W^{1,\infty}(\I,\R^d)$ such that $\dot{q}\in BV(\I,\R^d)$ such that
\be{eq:includiff} \left\{ \begin{array}{l} 
                  d\dot{q} + \NN(C(t),q)  \ni f(t,q) dt  \vsp \\
                  \forall t\in \overset{\circ}{\I},\ \dot{q}(t^+)=\PPP_{\C_{t,q(t)}}[\dot{q}(t^-)] \vsp \\
                  q(0)=q_0 , \vsp \\
                  \dot{q}(0)=u_0
                 \end{array} \right.\ee
where $q_0\in \textrm{Int}[C(0)]$ and $u_0$ are initial data and $\NN(C(t),q)$ is the proximal normal cone to $C(t)$ at $q$ (see Definition \ref{def:coneprox}).
Let us give a more precise sense to this differential problem.

\begin{df} \label{def:solution} Let $\I:=[0,T]$ be a bounded time-interval.
A continuous function $q:\I\rightarrow \R^d$ is a solution of (\ref{eq:includiff}) if there exists another function $k:\I\rightarrow \R^d$ such that~:
\begin{enumerate}
 \item[a)] $q$ belongs to $W^{1,\infty}(\I,\R^d)$
 \item[b)] $\dot{q}$ and $k$ belong to $BV(\I,\R^d)$
 \item[c)] the following differential equation is satisfied in the sense of time-measures
\be{eq2sto} d\dot{q} + dk = f(t,q)dt  \ee
 \item[d)] for all $t\in \overset{\circ}{\I}$, the impact law $\dot{q}(t^+)=\PPP_{\C_{t,q(t)}}[\dot{q}(t^-)]$ holds
 \item[e)] the initial conditions are verified : $q(0)=q_0$ and $\dot{q}(0)=u_0$
 \item[f)] the differential measure $dk$ is supported on $\left\{t,\ q(t)\in\partial C(t)\right\}$:
\be{eqK2}  |k|(t) = \int_0^t {\bf 1}_{q(s) \in \partial C(s)} d|k|(s), \qquad k(t)=\int_0^t \xi(s) d|k|(s),\ee
where $\xi:\I\rightarrow \R^d$ is a measurable function satisfying for all $s\in \I$: $\xi(s)\in \NN(C(s),q(s))$, $|\xi(s)|=1$ and $|k|(t):= \textrm{Var}\, (k,[0,t])$.
\end{enumerate}
\end{df}

\mb The main subject of this work is to prove the following global existence result.

\begin{thm} \label{thm:cv} Let $\I=[0,T]$ be a bounded time interval. Suppose that $C:\I \rightrightarrows \R^d$ is a Lipschitz and admissible (see Definition \ref{def:adm}) set-valued map
Assume that $f:\I \times \R^d \rightarrow \R^d $ is a measurable map satisfying:
\begin{align}
  & \exists K_L >0 \virg \forall t\in \I\virg \forall q,\tilde{q} \in C(t) \virg | f(t,q) -f(t,\tilde{q})| \leq K_L |q - \tilde{q}| &
\label{flip} \vsp \\
& \exists F \in L^1(\I) \virg \forall t\in \I\virg \forall q\in C(t) \virg |f(t,q)|\leq F(t).
\label{lingro}
\end{align}
Then the differential inclusion (\ref{eq:includiff}) admits at least one solution. 
\end{thm}


\mb The paper is structured as follows: the notions of uniform prox-regularity and admissibility related to the main assumptions of Theorem \ref{thm:cv} are detailed in Section \ref{sec:def}. Moreover, we prove a natural expression of the set $\C_{t_0,q_0}$ with the help of the tangent cone of $\Omega:=\{(t,q)\in \R \times \R^d,\ t\in \I,\ q\in C(t)\}$ at $(t_0,q_0)$. Then Section \ref{sec:global} is devoted to the proof of Theorem \ref{thm:cv}. In Section \ref{sec:kpart}, we deal with a set-valued map $C(\cdot)$, defined by constraints as in (\ref{C:interb}). More precisely, we prove under (\ref{reverse:intro}) global existence results (Theorem \ref{thm:particular}). For time-independent constraints, Assumption ($R_0$) is shown to be sufficient to obtain local existence results (see Theorem \ref{thm:partlocal}), which 
is weaker than the usual geometrical assumptions.

\section{Preliminaries and Definitions} \label{sec:def}

\mb First we precise some notations. For a time interval $\I$, we write $W^{1,\infty}(\I,\R^d)$ (resp. $W^{1,1}(\I,\R^d)$) for the Sobolev space of functions in $L^\infty(\I,\R^d)$ (resp. $L^1(\I,\R^d)$) whose derivative is also in $L^\infty(\I,\R^d)$ (resp. $L^1(\I,\R^d)$). $BV(\I,\R^d)$ is the space of functions in $L^\infty (\I,\R^d)$ with bounded variations on $\I$. We define the dual space ${\mathcal M} (\I) := ({\mathcal C}_c(\I))'$ where ${\mathcal C}_c(\I)$ is the space of continuous functions with compact support (corresponding to the set of Radon measure due to Riesz Theorem). We set ${\mathcal M}_+ (\I)$ for the subset of positive measures.

\mb
We emphasize that the different notions defined in this section can be extended in the case of an infinite dimensional Hilbert space $H$. We consider the Euclidean space $\R^d$, equipped with its euclidean metric $|\ |$ and its inner product $\langle\cdot,\cdot \rangle$. For a subset $Q$ of $\R^d$, we denote by $d_Q$ the distance function to this set~:
$$d_Q(x) := \inf_{y\in Q} \ |y-x|.$$ 

\begin{df} Given a family of sets $D_\epsilon\subset \R^d$ indexed by $\epsilon>0$, the outer and inner limits are defined respectively by
$$ \limsup_{\epsilon \searrow 0} D_\epsilon := \left\{ x\in \R^d,\ \exists x_k\to x,\ \exists \epsilon_k\searrow 0 \textrm{ with } x_k\in D_{\epsilon_k} \right\},$$
and
$$ \liminf_{\epsilon \searrow 0} D_\epsilon := \left\{ x=\lim_{\epsilon \to 0} x_\epsilon,\ x_\epsilon \in D_{\epsilon}\right\}.$$
If these two sets both equal a set $D\subset \R^d$, then we say that $D_\epsilon$ converges to $D$ and we write $\dsp \lim_{\epsilon \searrow 0} D_\epsilon=D$.
\end{df}

\begin{df} Let $Q$ be a closed subset of $\R^d$. The set-valued projection operator $P_Q$ is defined on $\R^d$ by
$$ \forall x\in \R^d, \qquad P_Q(x):=\left\{y\in Q,\ |x-y|=d_Q(x)\right\}.$$
\end{df}

\begin{df} \label{def:coneprox} Let $Q$ be a closed subset of $\R^d$ and $x\in Q$, we write $\NN(Q,x)$ for the proximal normal cone of $Q$ at $x$, defined by:
 $$ \NN(Q,x):=\left\{v\in\R^d,\ \exists s>0, \ x\in P_Q(x+sv)\right\}.$$
\end{df}

\begin{figure}
\psfragscanon
\psfrag{x0}[l]{$x_0$}
\psfrag{x1}[l]{$x_1$}
\psfrag{x2}[l]{$x_2$}
\psfrag{x3}[l]{$x_3$}
\psfrag{x4}[l]{$x_4$}
\psfrag{N0}[l]{$\NN(Q,x_0)$}
\psfrag{N1}[l]{$\NN(Q,x_1)$}
\psfrag{N3}[l]{$\NN(Q,x_3)$}
\psfrag{N4}[l]{$\NN(Q,x_4)$}
\psfrag{C}[l]{$Q$}
\begin{center}
\resizebox{0.5\textwidth}{!}{\includegraphics{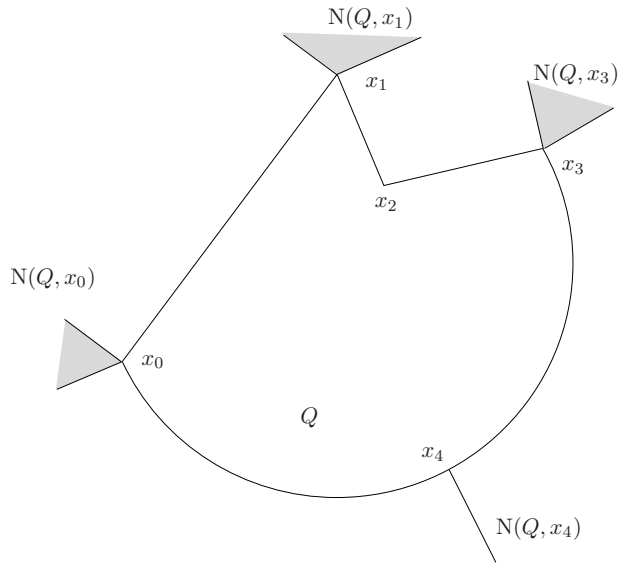}}
\end{center}
\caption{Examples of proximal normal cones}
\label{fig1}
\end{figure}

\mb The cone $\NN(Q,x)$ somehow generalizes the notion of the outward normal direction. In Figure \ref{fig1}, we have plotted a set $Q$ with several points $x_i$. At the regular point $x_4$ (where the boundary is smooth), the proximal normal cone is exactly the half-line directed by the outward direction. At the points $x_0,x_1$ and $x_3$ the boundary is not smooth and the proximal normal vectors constitute a cone. At the point $x_2$, notice that the proximal normal cone is reduced to $\{0\}$.

\mb We now define the tangent cone as follows:

\begin{df} Let $Q$ be a closed subset of $\R^d$ and $x\in Q$, we write $\TT_Q(x)$ for the tangential cone of $Q$ at $x$, defined by the following outer limit:
 $$ \TT_Q(x):=\limsup_{\epsilon \searrow 0} \, \frac{Q-x}{\epsilon} = \left\{v\in\R^d,\ \exists v_k \to v,\ \exists \epsilon_k \searrow 0, \  v_k \in \frac{Q-x}{\epsilon_k} \right\}.$$
\end{df}

\subsection{Uniform Prox-regularity}

\mb We now come to the main notion of ``prox-regularity''. It was initially introduced by H. Federer (in \cite{Federer}) in finite dimensional spaces under the name of ``positively reached sets''. Then it was extended in infinite dimensional space and studied by F.H.~Clarke, R.J.~Stern and P.R.~Wolenski in \cite{Clarke} and by R.A. Poliquin, R.T.~Rockafellar and L.~Thibault in \cite{PRT}.

\begin{df} Let $Q$ be a closed subset of $\R^d$ and $\eta>0$. The set $Q$ is said $\eta$-prox-regular if for all $x\in Q$ and $v\in \NN(Q,x)\setminus \{0\}$
 $$B\left(x+\eta\frac{v}{|v|},\eta\right) \cap Q = \emptyset.$$ By extension, a closed set $Q$ is said uniformly prox-regular if there exists $\eta>0$ such that $Q$ is $\eta$-prox-regular.
\end{df}

\mb We refer the reader to \cite{Clarke,Clarke2} for other equivalent definitions related to the limiting normal cone. The previous definition is very geometric, it describes the fact that we can continuously roll an external ball of radius $\eta$ on the whole boundary of the set $Q$. The main property is the following one: for a $\eta$-prox-regular set $Q$, and every $x$ satisfying $d_Q(x)<\eta$, the projection of $x$ onto $Q$ is well-defined (i.e. $P_Q(x)$ is a singleton) and the projection is continuous. \\
The notion of prox-regularity can be described by the hypomonotonicity property of the proximal normal cone too  (see the work of R.A. Poliquin, R.T. Rockafellar and L. Thibault (\cite{PRT})):

\begin{prop}[\cite{PRT}] \label{prop:hypo} Let $Q$ be a closed subset of $\R^d$. Then $Q$ is $\eta$-prox-regular if and only if the proximal normal cone $\NN(Q,\cdot)$ is an hypomonotone operator with constant $\frac{1}{2\eta}$. 
That means: for all $x,y\in Q$ and all $(\alpha,\beta)\in \NN(Q,x)\times \NN(Q,y)$, we have
$$ \langle \alpha-\beta,x-y\rangle \geq -\frac{1}{2\eta}\left[|\alpha|+|\beta|\right] |x-y|^2$$
or equivalently
\be{hypo}
\langle \alpha,y-x\rangle \leq \frac{1}{2\eta}|\alpha| |x-y|^2.
\ee
\end{prop}

\begin{prop} \label{prop:hypo2} Let $Q$ be a $\eta$-prox-regular set of $\R^d$. Consider $y\in \R^d$ with $d_{Q}(y)<2\eta$. Then $y \in Q$ if and only if
$$  \forall x\in Q\cap \overline{B}(y,d_Q(y)),\ \forall \alpha\in \NN(Q,x),\  \langle \alpha,y-x\rangle \leq \frac{1}{2\eta}|\alpha| |x-y|^2 .$$
\end{prop}

\dem Consider a point $y\notin Q$ satisfying the above property and choose $x\in \PPP_Q(y)$. Then $y-x\in \NN(Q,x)$ and we get
$$ |y-x|^2 \leq \frac{1}{2\eta}|x-y|^3,$$
which leads to a contradiction with the assumption $d_{Q}(y)<2\eta$.
\findem

\mb Since the prox-regularity of a set implies its Clarke's regularity (see Definition 6.4 of \cite{RW}), we have the following lemma (Corollary 6.29 of \cite{RW}).

\begin{lem} \label{lem:cone} Let $Q$ be a uniformly prox-regular set. Then for all $x\in Q$ the normal and tangential cones are mutually polar :
$$ \NN(Q,x)^\circ = \TT_Q(x).$$
That means:
$$ \TT_Q(x) =\left\{ u\in \R^d,\ \forall v\in \NN(Q,x) \ \langle u,v\rangle \leq 0 \right\}.$$
\end{lem}

\mb In particular, it comes that every tangential cone of a uniformly prox-regular set is convex. 

\begin{prop} \label{prop:Cconstant} Let $C(\cdot)$ be a Lipschitz continuous set-valued map taking uniformly prox-regular values on $\I$. For all $t_0\in\I$ and $q_0\in C(t_0)$, we set
$$ \C_{t_0,q_0}:=\liminf_{h\searrow 0} \ \frac{C(t_0+h)-q_0}{h} \quad \textrm{ and }\quad \tilde{\C}_{t_0,q_0} := \liminf_{\genfrac{}{}{0pt}{}{h\to 0}{\genfrac{}{}{0pt}{}{q \to q_0}{\genfrac{}{}{0pt}{}{t\to t_0}{q\in C(t)}}}} \ \frac{C(t+h)-q}{h}.$$
Then $\C_{t_0,q_0}=\tilde{\C}_{t_0,q_0}$.
\end{prop}

\dem Let $t_0\in\I$ and $q_0\in C(t_0)$, the first inclusion $\C_{t_0,q_0} \subset \tilde{\C}_{t_0,q_0}$ is obvious and we only deal with the other one. So let us fix $v\in \tilde{\C}_{t_0, q_0}$. By definition, there exist vectors $v_{h,t,q} \in \frac{C(t+h)-q}{h}$ such that
$$ v= \lim_{\genfrac{}{}{0pt}{}{h\to 0}{\genfrac{}{}{0pt}{}{t\to t_0}{q\to q_0}}} v_{h,t,q}.$$ We write
$$ z_{h,t,q} \in \PPP_{C(t_0+h)}(q_0+hv_{h,t,q}) \quad \textrm{and} \quad w_{h,t,q} := \frac{z_{h,t,q}-q_0}{h} \in \frac{C(t_0+h)-q_0}{h}.$$
It comes 
\be{eq:disbis} |v_{h,t,q}-w_{h,t,q}| = \frac{|hv_{h,t,q}-z_{h,t,q}+q_0|}{h} = \frac{d_{C(t_0+h)}(q_0+hv_{h,t,q})}{h}. \ee
Moreover $q_0+hv_{h,t,q}-z_{h,t,q} \in \NN(C(t_0+h),z_{h,t,q})$. Let $\xi_{h,t,q}\in \PPP_{C(t_0+h)}(q+hv_{h,t,q})$,
we have
\be{eq:hypp} |q_0+h v_{h,t,q} - z_{h,t,q}|\leq |\xi_{h,t,q} - z_{h,t,q}| + |\xi_{h,t,q}-(q_0+hv_{h,t,q})|. \ee
Then for $h>0$, let us choose $t_h\in\I$ such that $|t_h-t_0|\leq h^2$ and $q_h\in C(t_h)$ satisfying $|q_h-q_0|\leq h^2$. So by this way, 
\begin{align}
 |\xi_{h,t_h,q_h}-(q_0+hv_{h,t_h,q_h})| & \leq  |\xi_{h,t_h,q_h}-(q_h+hv_{h,t_h,q_h})| + |q_h-q_0| \nonumber \\
 & \leq d_{C(t_0+h)}(q_h+hv_{h,t_h,q_h}) +  |q_h-q_0| \nonumber \\
 & \leq d_H(C(t_0+h),C(t_h+h)) + |q_h-q_0| \nonumber \\
 & \leq c_0 |t_h-t_0| + |q_h-q_0| \leq (c_0+1) h^2. \label{j1}
\end{align}
In addition,
\begin{align*}
 d_{C(t_0+h)}(q_0+h v_{h,t_h,q_h}) & \leq d_{C(t_0)}(q_0+h v_{h,t_h,q_h}) + d_H(C(t_0+h),C(t_0)) \\
 & \leq h |v_{h,t_h,q_h}| + c_0 h
\end{align*}
and
\begin{align*}
 d_{C(t_0+h)}(q_h+h v_{h,t_h,q_h}) & \leq d_{C(t_h)}(q_h+h v_{h,t_h,q_h}) + d_H(C(t_0+h),C(t_h)) \nonumber \\
 & \leq h |v_{h,t_h,q_h}| + c_0 |t_h-t_0-h|. 
\end{align*}
For $h$ small enough, it comes 
$$ d_{C(t_0+h)}(q_0+h v_{h,t_h,q_h})\leq \frac{\eta}{2} \quad \textrm{and} \quad d_{C(t_0+h)}(q_h+h v_{h,t_h,q_h})\leq \frac{\eta}{2}.$$
By Theorem 4.8 in \cite{Clarke}, since $C(t_0+h)$ is $\eta$-prox-regular, the projection operator is Lipschitz on a neighborhood of $C(t_0+h)$. As a consequence, for $h$ small enough
\begin{align}
 |\xi_{h,t_h,q_h} - z_{h,q_h}| & = \left|\PPP_{C(t_0+h)}(q_h+h v_{h,t_h,q_h}) - \PPP_{C(t_0+h)}(q_0+h v_{h,t_h,q_h}) \right| \nonumber \\
& \leq \frac{2\eta}{2\eta-d_{C(t_0+h)}(q_0+h v_{h,t_h,q_h}) - d_{C(t_0+h)}(q_h+h v_{h,t_h,q_h})}|q_h-q_0| \nonumber \\
 & \leq 2|q_h-q_0| \leq 2 h^2. \label{j2}
\end{align}
It follows from (\ref{eq:hypp}), (\ref{j1}) and (\ref{j2})
\begin{align*}
 |q_0+h v_{h,t_h, q_h} - z_{h,q_h}| &  \leq |\xi_{h,t_h,q_h} - z_{h,t_h,q_h}| + |\xi_{h,t_h,q_h}-(q_0+hv_{h,t_h,q_h})| \\
 & \leq 2h^2 + (c_0+1) h^2 \leq C h^2,
\end{align*}
for some numerical constant $C$. It follows from (\ref{eq:disbis}) that 
$$ |v_{h,t_h,q_h}-w_{h,t_h,q_h}| \leq Ch$$
and so 
$$ v=\lim_{h\to 0} v_{h,t_h,q_h} = \lim_{h\to 0} w_{h,t_h,q_h}$$
with $w_{h,t_h,q_h} \in \frac{C(t_0+h)-q_0}{h}$. Hence $v\in \C_{t_0,q_0}$, which concludes the proof of this inclusion $\tilde{\C}_{t_0,q_0} \subset \C_{t_0,q_0}$. 
\findem


\subsection{The sets of admissible velocities}

In this subsection, we consider a set-valued map $C(\cdot)$ Lipschitz continuous and taking (nonempty) uniformly prox-regular values. We aim to describe the sets of admissible velocity $\C_{t,q}$ with the help of the whole set 
$$ \Omega:= \left\{ (t,q),\ t\in \I \ \textrm{and} \ q\in C(t) \right\}$$
and to prove its convexity.

\mb Let first recall the notion of  the derivable tangent cone (see Section 6.A in \cite{RW}).

\begin{df}
Let $Q$ be a closed set of $\R^d$. A vector $v\in \TT_Q(x)$ is said derivable if for all small enough $\epsilon >0$ there exists $v_\epsilon\in \frac{Q-x}{\epsilon} $ such that $v_\epsilon$ converges to $v$.
We note $\TT^D_Q(x)\subset \TT_Q(x)$ the set of derivable vectors, which can be seen as the following inner limit
$$ \TT^D_Q(x):= \liminf_{\epsilon \searrow 0} \, \frac{Q-x}{\epsilon} = \left\{v\in\R^d,\ v=\lim_{\epsilon \searrow 0} v_\epsilon \ \textrm{ with } \ v_\epsilon \in \frac{Q-x}{\epsilon} \right\}.$$
\end{df}

\mb A closed subset $Q$ is said {\it geometrically derivable} if every tangent vector is derivable, i.e. for all $x\in Q$
$$ \TT_Q(x)=\TT^D_Q(x).$$
In this particular case, the inner and outer limits are equal and so the limit is well-defined:
$$ \TT_Q(x)= \TT^D_Q(x) = \lim_{\epsilon \searrow 0} \, \frac{Q-x}{\epsilon}.$$

\mb For example, it is well-known that every uniformly prox-regular or even Clarke's regular set is geometrically derivable (see Corollary 6.30 in \cite{RW}).

\mb These definitions allow us to describe the set of admissible velocities with the the tangent cone of $\Omega$. For $(t,q)\in \Omega$, we have defined the set of admissible velocities
$$ \C_{t,q}:=\left\{v =\lim_{\epsilon \searrow 0} v_\epsilon, \ \textrm{ with } \ v_\epsilon \in \frac{C(t+\epsilon)-q}{\epsilon}\right\}.$$ 

\begin{prop} \label{prop:conebis} For every $(t,q)\in \Omega$ with $t\in \overset{\circ}{\I}$, we have
\be{eq:conebis}  \C_{t,q}=\left\{ u\in \R^d,\ (1,u) \in  \TT^D_\Omega((t,q))\right\} . \ee
\end{prop}

\dem Let us denote by $\widetilde{\C}_{t,q}$ the set of the right side in (\ref{eq:conebis}). It is obvious that for $v\in \C_{t,q}$, the vector $(1,v)$ is tangent to $\Omega$ at $(t,q)$ and derivable. Consequently the inclusion 
$\C_{t,q} \subset \widetilde{\C}_{t,q} $ is proved. \\
Let us now study the inverse inclusion. So consider $v\in \widetilde{\C}_{t,q}$. By definition $(1,v)$ is a derivable vector to $\Omega$ at $(t,q)$ so for small enough $\epsilon>0$, there exists a vector $(s_\epsilon,v_\epsilon)$ converging to $(1,v)$ such that
$$ (s_\epsilon,v_\epsilon)\in \frac{\Omega- (t,q)}{\epsilon}$$
which is equivalent to
$$ q+ \epsilon v_\epsilon \in C(t+\epsilon s_\epsilon).$$
Since the set-valued map $C$ is Lipschitz continuous (with a Lipschitz constant $c_0$), it comes 
\be{eq:dist} d_{C(t+\epsilon)}(q+\epsilon v_\epsilon) \leq d_H(C(t+\epsilon s_\epsilon), C(t+\epsilon)) \leq c_0 \epsilon |s_\epsilon-1|. \ee
So let us denote $z_\epsilon\in \PPP_{C(t+\epsilon)}(q+\epsilon v_\epsilon)$ and choose $w_\epsilon$ such that
$$ z_\epsilon := q+\epsilon w_\epsilon.$$
By definition, $z_\epsilon \in C(t+\epsilon)$, which means
$$ (1,w_\epsilon) \in \frac{\Omega- (t,q)}{\epsilon}.$$
Moreover (\ref{eq:dist}) yields
$$ |v_\epsilon-w_\epsilon| \leq c_0|s_\epsilon -1|.$$
Since $s_\epsilon$ converges to $1$ and $v_\epsilon$ to $v$, we deduce that $w_\epsilon$ converges to $v$ too. Thus $v\in \C_{t,q}$, which ends the proof of the inclusion $ \widetilde{\C}_{t,q} \subset \C_{t,q}$.
\findem

\mb Then we can state the main result concerning the sets of admissible velocities.

\begin{prop}  \label{prop:coneconvex} For every $(t,q)\in \Omega$ with $t\in \overset{\circ}{\I}$, the set $\C_{t,q}$ is convex.
\end{prop}

\dem If $\C_{t,q}$ is supposed to be nonempty, we can choose two vectors $v^1$ and $v^2$ belonging to $\C_{t,q}$. Let $\alpha \in ]0,1[$, we aim to prove that
$$ v:= \alpha v^1+ (1-\alpha) v^2\in \C_{t,q}.$$
According to the previous lemma, for $i=1,2$ and small enough parameter $\epsilon$, we have a sequence $v^i_\epsilon$, converging to $v^i$ such that
$$  q+\epsilon v^i_\epsilon \in  C(t+\epsilon).$$
Writing $v_\epsilon:= \alpha v^1_\epsilon + (1-\alpha) v^2_\epsilon$, we want to estimate $d_{C(t+\epsilon)}(q+\epsilon v_\epsilon)$. Suppose that $q+\epsilon v_\epsilon$ does not belong to $C(t+\epsilon)$ and let $z_\epsilon \in \PPP_{C(t+\epsilon)}(q+\epsilon v_\epsilon)$. So $q+\epsilon v_\epsilon-z_\epsilon$ is a proximal normal vector at $z_\epsilon$. Since $C(t+\epsilon)$ is $\eta$-prox-regular, thanks to the hypomonotonicity property of the proximal normal cone (see Proposition \ref{prop:hypo}), we have for $i=1,2$
$$ \langle q+\epsilon v^i_\epsilon - z_\epsilon, q+\epsilon v_\epsilon-z_\epsilon \rangle \leq \frac{|q+\epsilon v_\epsilon-z_\epsilon|}{2\eta}|q+\epsilon v^i_\epsilon - z_\epsilon|^2.$$
Multiplying by $\alpha$ the previous inequality for $i=1$ and by $1-\alpha$ the one for $i=2$ and summing them, we obtain
$$ |q+\epsilon v_\epsilon-z_\epsilon| \leq \frac{1}{2\eta}\left[ \alpha |q+\epsilon v^1_\epsilon - z_\epsilon|^2+ (1-\alpha) |q+\epsilon v^2_\epsilon - z_\epsilon|^2\right].$$
Since $q+\epsilon v^i_\epsilon \in C(t+\epsilon)$ for $i=1,2$, we have
\begin{align*}
  |q+\epsilon v^i_\epsilon - z_\epsilon| & \leq |q+\epsilon v_\epsilon -z_\epsilon| + \epsilon |v_\epsilon - v^i_\epsilon| \\
   &  \leq 2\epsilon |v_\epsilon - v^i_\epsilon| \\
   & \leq 2\epsilon (|v^1_\epsilon |+ |v^2_\epsilon |) \leq 4M\epsilon,
\end{align*}
where $M$ denotes a uniform bound of $v^1_\epsilon$ and $v^2_\epsilon$ (since they are convergent). So we deduce that
\be{eq:dist2} |q+\epsilon v_\epsilon-z_\epsilon| \leq \frac{8M^2}{\eta} \epsilon^ 2.\ee
Then let us choose $w_\epsilon$ such that 
$$ z_\epsilon = q+\epsilon w_\epsilon \in C(t+\epsilon).$$
From (\ref{eq:dist2}), we know that
$$ \epsilon |v_\epsilon-w_\epsilon | \leq \frac{8M^2}{\eta} \epsilon^ 2$$
and so $w_\epsilon$ converges to $v$. We have proved that $v=\alpha v^1+ (1-\alpha) v^2$ is the limit of velocities $w_\epsilon$ satisfying
$$ q+\epsilon w_\epsilon \in C(t+\epsilon).$$
That means $v\in \C_{t,q}$, which shows the convexity of this cone. \findem

\subsection{Admissibility}

After having introduced the notion of prox-regularity, let us now present the concept of admissibility.

\begin{df} \label{def:adm} A set-valued map $C:\I\rightrightarrows \R^d$ is said {\em admissible} on $\I:=[0,T]$ if it takes uniformly prox-regular values (with a same constant) and if there exist $\delta,r,\tau>0$, and for all $t_0\in[0,T]$ sequences $(x_p)_p$ and $(u_p)_p$ with $|u_p|=1$ and $x_p\in C(t_0)$ such that for all $s\in [0,T]$ with $|t_0-s|\leq \tau$, $(B(x_p,r))_p$ is a bounded covering of the boundary $\partial C(s)$ and 
\be{hyp5} \forall p,\ \forall x\in \partial C(s) \cap B(x_p,2r),\ \forall v\in \NN(C(s),x), \quad \langle v,u_p\rangle \geq \delta |v|.\ee
\end{df}

\mb For an admissible set, the ``good directions'' $u_p$ allow us to build inward cones.

\begin{lem} Let us consider an admissible Lipschitz set-valued map $C:\I\rightrightarrows \R^d$ and let us keep the notations of the previous definition. \\
Fix $t_0\in \I$, $x_p\in C(t_0)$ and $\nu>0$. Then if $t\in\I$ with $|t_0-t|+\nu\leq \tau$ and $x \in B(x_p,3r/2) \bigcap C(t)$, $$ \overline{B}\left(x-\nu \kappa_0 u_p,\frac{\nu\delta}{2} \right) \subset C(t+\nu),$$
 as soon as $ \nu <  \nu_{min}:=\dsp \min \left\{\dfrac{\eta \delta}{ (2\kappa_0+2c_0+\delta)^2}\ ,\ \dfrac{r}{2(c_0+\delta+2\kappa_0)} \right\}$ and $\kappa_0:=\frac{c_0}{\delta}+1$,
 where $c_0$ is the Lipschitz constant of the set-valued map $C(\cdot)$ (see (\ref{eq:Qlip})) and $\eta$ the prox-regularity constant.
\label{lem:conerentrant}
\end{lem}

\dem Let $\nu < \nu_{min}$ and $x \in B(x_p,3r/2) \bigcap C(t) $. We consider $z=x-\nu \kappa_0 u_p+\theta$, with $|\theta | \leq \frac{\nu\delta}{2}$.
Suppose that $z\notin C(t+ \nu)$ and set $y \in \PPP_{C(t+\nu)}(z)$, $v:= z-y \in \NN(C(t+\nu),y) \setminus \{0\}$. 
Necessarily, $y $ belongs to $B(x,r/2)$ because 
$$|x-y|\leq |x-z|+|z-y|\leq 2|x-z|+c_0 \nu\leq  2\left(\kappa_0\nu + \dfrac{\nu\delta}{2} \right) +c_0\nu \leq\nu(\delta+2\kappa_0 +c_0)< r/2 ,$$
where we have used that $d_{C(t+\nu)}(z)\leq d_{C(t)}(z)+d_H(C(t),C(t+\nu))$ and $x\in C(t)$. \\
Let$\xi\in \PPP_{C(t+ \nu)}(x)$, since $C(t+\nu)$ is $\eta$-prox-regular, it comes 
 $$\ww  \xi -y,v\xx  \leq \frac{1}{2\eta}|\xi-y|^2|v|.$$ 
That yields 
\begin{align*}
\ww  \xi -z,v\xx + |v|^2  & \leq \frac{1}{2\eta}\left(|\xi -x|+ |x-y| \right)^2|v| \\
& \leq \frac{1}{2\eta}\nu^2(\delta+2\kappa_0+ 2c_0)^2|v|. 
\end{align*}
Consequently, 
\begin{align} 
\lefteqn{\kappa_0 \ww  \nu u_p, v\xx - \ww \theta, v\xx } & & \nonumber \\
& & \leq \frac{1}{2\eta}(\nu(\delta+2\kappa_0+ 2c_0))^2|v|+d_{C(t+\nu)}(x)|v| \leq \left[ \frac{1}{2\eta}(\nu(\delta+2\kappa_0+2c_0))^2 + c_0 \nu \right] |v| . \label{ineg1}
\end{align}
As $y \in B(x,r/2) $, $y $ belongs to $B(x_p,2r)$ and the admissibility assumption gives $\ww  u_p, v\xx \geq \delta |v|$.
Thus (\ref{ineg1}) implies that $$ \kappa_0 \delta \nu  - \nu \frac{\delta}{2} \leq \frac{1}{2\eta}(\nu(\delta+2\kappa_0+ 2c_0 ))^2 +c_0 \nu  $$
and so that $\nu \geq  \frac{\eta \delta}{ (2\kappa_0+\delta+ 2c_0)^2}$, which leads to a contradiction. Thus $z\in C(t+\nu)$.
\findem

\mb The following lemma is a consequence of the previous one. 

\begin{lem} \label{lem:conerentrant2} Under the previous assumptions and notations, fix $t_0\in \I$ and $x_p\in C(t_0)$.
For all $t\in \I$ and $x \in B(x_p,3r/2) \bigcap C(t)$, for all $  \nu < \nu_{min}$,  $- \kappa_0 u_p \in \dfrac{C(t+ \nu) -x}{\nu} $. \newline
Moreover for all $ y \in \frac{C(t+\nu) -x}{\nu} \bigcap B(0,\frac{r}{2\nu})$,  
$$-\kappa_0 u_p \in \TT_{\frac{C(t+\nu)-x}{\nu}} (y) . $$
In addition, $-\kappa_0 u_p-y$ belongs to $\TT_{\frac{C(t+\nu) -x}{\nu}}(y)$ as soon as
\be{eq:as} \nu |y|^2 \leq \delta \eta/2. \ee
\end{lem}

\dem 
Concerning the first point, by Lemma \ref{lem:conerentrant}, $ x -\kappa_0 \nu u_p \in C(t+\nu)$ if $ x\in B(x_p,3r/2) \bigcap C(t)$ and $\nu <  \nu_{min}$. Thus $-\kappa_0 u_p \in \dfrac{C(t+\nu) -x}{\nu} $. 
Let $ y \in \frac{C(t+\nu) -x}{\nu} \bigcap B(0,\frac{r}{2\nu})$,  $x +\nu y \in B(x,r/2) \subset B(x_p,2 r)  $ and $x +\nu y \in C(t+\nu) $. 
The admissibility assumption implies that for all $v \in \NN(C(t+\nu) ,x+\nu y)$, $$\ww -u_p,v \xx \leq -\delta |v| \leq 0. $$
Since $\NN(C(t+\nu),x+\nu y)= \NN( \frac{C(t+\nu) -x}{\nu},y) $, we have for all $v \in \NN( \frac{C(t+\nu) -x}{\nu},y)$,$$\ww -u_p,v \xx \leq 0. $$
As a consequence, $-u_p \in  \NN( \frac{C(t+\nu)-x}{\nu},y)^\circ = \TT_{\frac{C(t+\nu)-x}{\nu}} (y) $ (due to Lemma \ref{lem:cone}). \\
Let us now prove the last point.  Since $\frac{C(t+\nu) -x}{\nu}$ is $\frac{\eta}{\nu}$-prox-regular and contains $a\in\frac{\PPP_{C(t+\nu)}(x)-x }{\nu}$ (satisfying $|a|\leq c_0$),  for all $y\in C(t+\nu)$ and $v\in \NN( \frac{C(t+\nu)-x}{\nu},y) = \NN(C(t+\nu) ,x+\nu y)$
$$ \langle a-y,v\rangle \leq \frac{\nu}{2\eta}|v||a-y|^2\leq  \frac{\nu}{\eta}|v| \left(c_0^2 + |y|^2\right) , $$
by the hypomonotonicity property of the proximal normal cone. It follows from (\ref{eq:as}) that
\begin{align*}
 \langle -\kappa_0 u_p-y,v\rangle & \leq  -\kappa_0 \delta |v| + \langle -y,v\rangle \leq -\kappa_0\delta |v| + \frac{\nu}{\eta}|v|\left(c_0^2 + |y|^2\right) + c_0 |v| \\
 & \leq -\frac{3\delta}{4}|v| +  \frac{\nu}{\eta}|v| |y|^2 \leq 0
\end{align*}
(since $\nu c_0^2\leq \eta \delta/4$ and $\kappa_0\delta=\delta+c_0$),  which proves the expected result (thanks to Lemma \ref{lem:cone}).
\findem

\begin{prop} \label{prop:assumption} Let $C(\cdot)$ be an admissible and Lipschitz continuous set-valued map on $\I$. For all $t_0\in\I$ and $q_0\in C(t_0)$, we set
$$ \C_{t_0,q_0}:=\liminf_{h\searrow 0} \ \frac{C(t_0+h)-q_0}{h} \quad \textrm{ and }\quad \tilde{\C}_{t_0,q_0} := \liminf_{\genfrac{}{}{0pt}{}{h\to 0}{\genfrac{}{}{0pt}{}{q \to q_0}{\genfrac{}{}{0pt}{}{t\to t_0}{q\in C(t)}}}} \ \frac{C(t+h)-q}{h}.$$
Then $\C_{t_0,q_0}=\tilde{\C}_{t_0,q_0}\neq \emptyset$.
\end{prop}

\dem We refer the reader to Proposition \ref{prop:Cconstant} for the equality $\C_{t_0,q_0}=\tilde{\C}_{t_0,q_0}$. It remains us to check that these sets are nonempty. \\
If $q_0\in \textrm{Int}(C(t_0))$ then we easily have $\C_{t_0,q_0}=\tilde{\C}_{t_0,q_0}=\R^d$. Else $q_0\in \partial C(t_0)$ and by Lemma \ref{lem:conerentrant}, there exist $u_p\neq 0$, $\kappa_0>0$ and $\nu_{min}>0$ such that for all $\nu\in]0,\nu_{min}[$, $q_0-\nu \kappa_0 u_p \in C(t_0+\nu)$. So we deduce that $-\kappa_0 u_p \in \C_{t_0,q_0}\neq \emptyset$.
\findem

\section{Discretization and convergence of approximate solutions} \label{sec:global}

This section is devoted to the proof of Theorem \ref{thm:cv}. As usual, we obtain existence results for (\ref{eq:includiff}) by proving the convergence of a sequence of discretized solutions. \\
To do so, we extend the {\it Catching-up algorithm} (proposed by J.J. Moreau for the firt-order differential inclusions) to the considered second-order problem. Let us describe the numerical scheme.

\medskip
Let $h:=T/N<\nu_{min}/4$ be the time step, where $\nu_{min}$ is defined in Lemma \ref{lem:conerentrant}. We denote by  $q_h^n\in\R^{d}$ and $u_h^n\in\R^{d}$ the approximated solution and velocity at time $t_h^n=nh$ for $n\in\{0,..,N\}$.

\mb The approximated solutions are built using the following scheme:
\begin{enumerate}
\item Initialization : \be{CondInit} (q_h^0,q_h^1):=(q_0,q_0+hu_0+h^2 f_h^0)\ \textrm{ with } \ f_h^0:=\frac{1}{h}\int_{0}^{t_h^{1}}f(s,q_h^0)ds.  \ee 
Since $q_0\in Int(C(0))$, $q_h^1$ belongs to $C(h)$ for $h$ small enough (such that $(|u_0|+\int_\I F(t) dt)h<d_{\partial C(0)}(q_0)$) .
\item Time iterations: $q_h^i$  are given for $i\in\{0, ... , n\}$.
We define $\dsp f_h^n:=\frac{1}{h}\int_{t_h^n}^{t_h^{n+1}}f(s,q_h^n)ds$ and 
\be{eq:qn+1} q_h^{n+1} \in \PPP_{C(t_h^{n+1})} \left[2q_h^n-q_h^{n-1} +h^2f_h^n\right]. \ee 
\end{enumerate}
This algorithm is a ``prediction-correction algorithm'': the predicted point $2q_h^n - q_h^{n-1}+h^2f_h^n$, that may not be admissible at the time $t_h^{n+1}$, is projected onto $C(t_h^{n+1})$. \\
We define the following functions: for all $t\in [t_h^n, t_h^{n+1}[$,
\be{eq:position} q_h(t) := q_h^n + (t-t_h^n) \frac{q_h^{n+1}-q_h^{n}}{h} \ee
and the velocity $u_h$ by
\be{eq:vitesse} u_h(t):=u_h^{n+1}:= \frac{q_h^{n+1}-q_h^n}{h}. \ee

 \mb Note that, for every $h<\frac{d_{\partial C(0)}(q_0)}{|u_0|+\int_\I F(t) dt}$, the scheme is well-defined. Moreover the computed configurations are feasible~:
\be{prop:feasible_config}
\forall h>0, \forall n\in\{0,...,N\},\quad  q_h(t_h^n)=q_h^n \in C(t_h^n).
\ee
 
\mb For the intermediate times $t\in ]t_h^n,t_h^{n+1}[$, the point $q_h(t)$ may not belong to $C(t)$. However from (\ref{prop:feasible_config}) and the Lipschitz regularity of the set-valued map $C(\cdot)$ (see (\ref{eq:Qlip})), we have the following estimate~:
\be{eq:feasible_config2} \forall h>0, \forall t\in \I,\quad d_{C(t)}(q_h(t)) \leq \max\{ d_{C(t)}(q_h^{n}),d_{C(t)}(q_h^{n+1})\} \leq c_0 h. \ee 

\gb The proof of Theorem \ref{thm:cv} is quite technical and will be decomposed into 4 steps, which we shall briefly describe below. 
\begin{itemize}
 \item In Subsection \ref{subsec:31}, we obtain uniform bounds on the computed velocities $u_h$ when $h$ goes to $0$ in $L^\infty(\I,\R^d)$ and in $BV(\I,\R^d)$.
\item In Subsection \ref{subsec:32}, we use compactness arguments in order to extract a subsequence of $(q_h,u_h)_{h>0}$ converging to $(q,u)$.
\item In Subsection \ref{subsec:33}, we check that the limit functions $q,u$ satisfy the second order differential equation (the momentum balance) and the initial conditions.
\item In Subsection \ref{subsec:34}, we verify the impact law for $u$.
\end{itemize}
After that, we will have proved that the limit function $q$ is a solution of Problem (\ref{eq:includiff}), which shows Theorem \ref{thm:cv}. \findem

\subsection{Uniform estimates on the computed velocities} \label{subsec:31}
This subsection is devoted to the proof of uniform estimates in $BV(\I,\R^d)$ for the computed velocities in order to extract a convergent subsequence by compactness arguments.


\mb First, we check that the velocities are uniformly bounded in $L^\infty(\I,\R^d)$. Aiming that, we prove the following lemma giving a first estimate on the velocities.

\begin{lem} \label{lem:velocity} For all integer $n\in\{0,... ,T/h\}$
$$ |u_h^{n+1}| \leq 2 |u_h^n+hf_h^n| + c_0. $$ 
\end{lem}

\dem By rewriting (\ref{eq:qn+1}) in terms of velocity, we deduce that
$$ u_h^{n+1} \in \PPP_{\frac{C(t_h^{n+1})-q_h^n}{h}}\left[ u_h^n+ hf_h^n\right].$$
With $z:=\PPP_{C(t_h^{n+1})}(q_h^n) $ (as $c_0h\leq c_0 \nu_{min}/4 \leq \eta$, the projection de $q_h^n$ is single-valued), it comes
$$ | u_h^n+ hf_h^n - u_h^{n+1} | \leq \left|u_h^n+ hf_h^n-\frac{z-q_h^n}{h}\right|.$$
The proof is also achieved thanks to $|z-q_h^n|\leq d_H(C(t_h^{n+1}),C(t_h^n)) \leq c_0 h$.
\findem

\mb To iterate this reasoning, it would be very interesting to obtain a similar estimate without the factor $2$ in Lemma \ref{lem:velocity}. This is the main difficulty in order to obtain a uniform bound of the velocities and we solve it in the following proof.

\begin{prop} \label{prop:uLI}
 There exists $h_1>0$ such that the sequence of computed velocities $(u_h)_{h<h_1}$ is bounded in $L^\infty(\I,\R^d)$. We set 
$$ K:= \sup_{h<h_1} \|u_h\|_{L^\infty(\I)} <\infty.$$ 
\end{prop}

\dem
\mb {\bf $1-)$} Estimate on the velocities for small time intervals.\\
Let us fix $h<h_1$ with $h_1$ later defined and satisfying $h_1\leq\min\{\nu_{min},\tau\}$. Consider a small time interval $[t_-,t_+]\subset I$ of length satisfying
\be{truc} h\leq |t_+-t_-|\leq \min\left\{ \frac{r}{5\left(|u_h^{n_0}|+2\kappa_0 + \int_0^T F(t) dt\right)}, \frac{\tau}{2} \right\} \ee
where $\kappa_0$ is introduced in Lemma \ref{lem:conerentrant} and $n_0$ is the smallest integer $n$ such that $t_h^n\geq t_-$. Thus $t_h^{n_0}\in [t_-,t_+]$.
We are looking for a bound on the velocity on this time interval. Rewriting the scheme in terms of velocities, we have that
$$ u_h^{n+1} \in \PPP_{\frac{C(t_h^{n+1})-q_h^n}{h}}\left[ u_h^n+ hf_h^n\right],$$
involving $u_h^n+ hf_h^n-u_h^{n+1}\in \NN\left(\frac{C(t_h^{n+1})-q_h^n}{h},u_h^{n+1}\right)$. \\
By the admissibility of $C$, there are balls $B(x_p,r)$, which cover $C(t_h^{n_0})$. So there exists at least one index $p$ such that $q_h^{n_0}\in B(x_p,r)$ and we denote $\omega:=u_p$. \\
Then for $n$ satisfying $t_h^n \in[t_-,t_+]$ (since $h+|t_h^n-t_h^{n_0}|\leq \tau$), Lemma \ref{lem:conerentrant2} yields that for $h< \nu_{min}$ if
\be{eq:error3} q_h^n\in B(q_h^{n_0},r/2)   \ee
(which implies $q_h^n\in B(x_p,3r/2)$) and
\be{eq:error3bis} h|u_h^{n+1}|^2\leq \eta\delta/2,\ee then
$ -\kappa_0 w-u_h^{n+1}\in \NN\left( \frac{C(t_h^{n+1}) -q_h^n}{h},u_h^{n+1} \right)^\circ$. So we deduce that
$$\langle  u_h^n+ hf_h^n - u_h^{n+1},-\kappa_0 w-u_h^{n+1}\rangle \leq 0,$$
which implies
$$|u_h^{n+1}+\kappa_0 w| \leq |u_h^n+ hf_h^n +\kappa_0 w|\leq |u_h^n +\kappa_0 w| + h|f_h^n|.$$
We set $m$ the smallest integer (bigger than $n_0$) such that $m+1$ does not satisfy (\ref{eq:error3}), (\ref{eq:error3bis}) or $t_h^{m+1}\notin [t_-,t_+]$. By summing these previous inequalities from $n=n_0$ to $n=p$ with $n_0\leq p\leq m$, we get
$$\forall p\in\{n_0,...,m\},\qquad |u_h^{p+1}+\kappa_0 w|\leq |u_h^{n_0}+\kappa_0 w| + \int_0^{t_h^{m+1}} F(t) dt.$$
Finally, it comes
\be{eq:error4} \sup_{n_0\leq p\leq m+1} |u_h^{p}| \leq |u_h^{n_0}|+2\kappa_0 + \int_0^{t_h^{m+1}} F(t) dt.\ee
By integrating in time, we deduce
$$ |q_h^{m+1}-q_h^{n_0}| \leq \left(|u_h^{n_0}|+2\kappa_0 + \int_0^{t_h^{m+1}} F(t) dt\right)(|t_+-t_-|+h)\leq \frac{2r}{5}<\frac{r}{2}, $$
by (\ref{truc}).
Consequently by definition of $m$, $t_h^{m}\leq t_+ < t_h^{m+1}$ as soon as (\ref{eq:error3bis}) holds for $n=m+1$. 
Moreover thanks to Lemma \ref{lem:velocity}, we have
$$ |u_h^{m+2}| \leq 2|u_h^{m+1}|+c_0+2\int_0^T F(t)dt.$$
From (\ref{eq:error4}), (\ref{eq:error3bis}) is satisfied for $n=m+1$ as soon as
$$ h \left(2(|u_h^{n_0}|+2\kappa_0 + \int_0^{T} F(t) dt)+c_0+2\int_0^T F(t) dt \right)^2 \leq \frac{\eta \delta}{2}.$$
Finally, we obtain
$$ \sup_{t_-\leq t_h^n\leq t_+} |u_h^{n}| \leq |u_h^{n_0}|+2\kappa_0 + \int_0^T F(t) dt,$$
as soon as 
$$ h \left(2(|u_h^{n_0}|+2\kappa_0 + \int_0^{T} F(t) dt)+c_0+2\int_0^T F(t) dt \right)^2 \leq \frac{\eta \delta}{2}.$$

\mb {\bf $2-)$} End of the proof.\\
Let $h<h_1$ (later defined), we are now looking for a bound on $u_h$ on the whole time interval $\I=[0,T]$, in assuming that $T<\tau/2$ without loss of generality. Let us start with $t_-=t(0):=0$. We know that with 
$$ t_+=t(1):= \min\left\{\frac{r}{5 A(1)},T \right\}$$
where 
$$A(1):=|u_0|+2\kappa_0 + \int_0^{T} F(t) dt$$
we have
$$ \sup_{0\leq t_h^n\leq t(1)} |u_h^n| \leq A(1) \leq |u_0|+2\kappa_0 + \int_0^T F(t) dt,$$
as soon as
$$  h \left( 2 A(1)+c_0+2\int_0^T F(t) dt \right)^2 \leq \frac{\eta \delta}{2}.$$
Then, let us suppose that there exists $n_1$ such that $t(0)<t^{n_1}_h\leq t(1)<t^{n_1+1}_h$. We have $0\leq\delta_1:=t(1)-t^{n_1}_h<h$. In that case, we set $t_-=t_h^{n_1}$ and  
\begin{align*} 
t_+=t(2)&:= \min\left\{t_h^{n_1}+\frac{r}{5 A(2)},T\right\}\\
&=\min\left\{t(1)-\delta_1+\frac{r}{5A(2)},T\right\},
\end{align*}
with 
$$ A(2) :=|u_0|+4\kappa_0  + 2\int_0^T F(t) dt \geq |u_h(t_h^{n_1})|+2\kappa_0  + \int_0^T F(t) dt.$$
From the previous point, we deduce that
\begin{align*}
 \sup_{t(1)\leq t_h^n\leq t(2)} |u_h^n| & \leq \sup_{t_h^{n_1}\leq t_h^n\leq t(2)} |u_h^n| \\ 
 & \leq |u_h^{n_1}|+2\kappa_0 + \int_0^T F(t) dt \\
& \leq |u_0| + 4\kappa_0  + 2\int_0^T F(t) dt=A(2)
\end{align*}
as soon as 
$$ h \left( 2 A(2)+c_0+2\int_0^T F(t) dt \right)^2 \leq \frac{\eta \delta}{2}.$$
Hence
$$ \sup_{0 \leq t_h^n\leq t(2)} |u_h^n| \leq A(2)=|u_0|+4\kappa_0 + 2\int_0^T F(t) dt.$$
By iterating this reasoning, for any integer $k\geq 1$ we set
$$ A(k):= |u_0|+2k \kappa_0 + k \int_0^{T} F(t)dt = A(k-1)+ 2\kappa_0 + \int_0^T F(t) dt$$
and
\begin{align*} 
t(k) & := \min\left\{t(k-1)-\delta_{k-1}+ \frac{r}{5A (k)},T \right\} \\
 & = \min\left\{-\sum_{i=1}^{k-1}\delta_i+\sum_{i=1}^k \frac{r}{5 A(i)},T\right\}.
\end{align*}
where $\delta_k<h$ for all $k$. This construction of $t(k)$ can be made while $t(k-1)-t(k-2)>h$. This condition will be verified as long as
$$ -\delta_{k-2}+ \frac{r}{5 A(k-1) }>h.$$
Using the fact that $0\leq \delta_{k-2}<h$, we see that we can construct $t(k)$ for $k<N$ verifying
$$\frac{r}{5\left(|u_0|+2(k-1)\kappa_0 + (k-1)\int_0^T F(t) dt\right)}>2h,$$
which is implied by
\be{k0}
k\leq k_0(h):=\left\lfloor 1+\left(\frac{r}{10h}-|u_0|\right) \left(2\kappa_0 + \int_0^T F(t) dt \right)^{-1} \right\rfloor.
\ee
Consequently, we know that the velocities can be bounded on $[0,t(k_0(h))]$ as follows
$$ \sup_{0\leq t_h^n \leq t(k_0(h))} |u_h^k| \leq A(k_0(h))$$
where
\be{eq:tk0}
t(k_0(h))=\min\left\{-\sum_{i=1}^{k_0(h)-1}\delta_i+\sum_{i=1}^{k_0(h)} \frac{r}{5 A(i)},T\right\},
\ee
under the property 
$$ h \left(2 \sup_{0\leq t_h^n\leq t(k_0(h))} |u_h^{n}| +c_0 +2\int_0^T F(t) dt\right) ^2 \leq h\left(2 A(k_0(h)) + c_0 + 2\int_0^T F(t) dt\right)^2 \leq \eta\delta/2. $$
Now, using the fact that $k_0(h)$ goes to infinity when $h$ goes to zero, that the harmonic serie diverges and that (\ref{k0}) yields
$$
\left|\sum_{i=1}^{k_0(h)-1}\delta_i\right|\leq h k_0(h)\leq C,
$$
we see by (\ref{eq:tk0}) that $t(k_0(h))$ is equal to $T$ for $h$ small enough. Therefore, there exists $h_0<\nu_{min}$ such that, for $h<h_0$, $T=t(k_0(h_0))=t(k_0(h))$. Finally, we see that, for $h<h_0$, $t(k)$ can be constructed until $k=k_0(h_0)$ (since $h\to k_0(h)$ is non decreasing) and that $t(k_0(h_0))=T$ (independently from $h<h_0$). Hence $u_h$ can be bounded as follows
\be{eq:final}  \sup_{0\leq t_h^n\leq T} \ |u_h^n| \leq A(k_0(h_0)), \ee
if
\be{eq:ass} h \left(2 A(k_0(h_0))+c_0 +2\int_0^T F(t) dt \right) ^2 \leq \eta\delta/2. \ee
However, (\ref{eq:ass}) holds for 
$$ h<h_1:=\min\left\{ h_0, \frac{\eta\delta}{ 2\left( 2 A(k_0(h_0))+c_0 +2\int_0^T F(t) dt  \right)^2} \right\}.$$
So we finally obtain the uniform bound
$$ \sup_{h<h_1} \ \sup_{0\leq t_h^n\leq T} \ |u_h^n| \leq A(k_0(h_0))= |u_0|+2k_0(h_0) \kappa_0 + k_0(h_0)\int_0^T F(t) dt, $$
which concludes the proof of a uniform bound in $L^\infty$ for the velocities $u_h$.
\findem

\mb Having obtain a uniform bound of the velocities, we can now prove that they have a uniformly bounded variation on the whole time-interval $I$.

\begin{thm}
\label{thm:uBV} There exists $h_2\in]0,h_1[$ such that the sequence of computed velocities $(u_h)_{h<h_2}$ is bounded in $BV(I,\R^d)$.
\end{thm}

\dem We adapt the proof of Theorem 3.2 in \cite{F-Aline}. To study the variation of $u_h$ on $I$, we split $I$ into smallest intervals. We define $(s_j)_j$ for $j$ from $0$ to $P$ such that:
$$
\left|\begin{array}{l}
\dsp s_0=0\virg s_P=T,\smallskip\\
\dsp|s_{j+1}-s_j|=\frac{1}{2}\min\left\{\tau, \frac{r}{K}\right\}, \hbox{ for } j=0\ldots P-2,\smallskip\\
\dsp|s_{P}-s_{P-1}|\leq \frac{1}{2}\min\left\{\tau, \frac{r}{K}\right\},
\end{array}\right.$$
where $\tau$ and $r$ are given by Definition~\ref{def:adm} of the admissibility and $K$
is the Lipschitz constant of $q_h$ (see Proposition
\ref{prop:uLI}). All these constants are  independent on $h$ and such a construction gives
\be{eq:boundP}
P=\left\lfloor \frac{2T}{\min\left\{\tau, \frac{r}{K}\right\}}\right\rfloor +1,
\ee 
which is independent of $h$.
Then, for all $h<h_1$, we define $n^j_h$ for $j$ from $0$ to $P-1$ as the first time step strictly greater than $s_j$: $$t_h^{n_h^j-1}\leq s_j < t_h^{n_h^j}, $$ and $n_h^P$ is set equal to $N$ ($t_h^N=t_h^{n_h^P}=T$). 

\mb
In the sequel, we suppose $h<\min\{|s_{j+1}-s_j|\}/2$. We also obtain a strictly increasing sequence of $(t_h^{n_h^j})_j$ verifying
\be{eq:dtnj}
|t_h^{n_h^j}-t_h^{n_h^{j-1}}|\leq \min\left\{\tau, \frac{r}{K}\right\}.
\ee

\mb The variation of $u_h$ on $I$ can be estimated as follows
$$
\textrm{Var}_I (u_h)=\sum_{n=0}^{N-1} |u_h^{n+1}-u_h^n| = \sum_{j=0}^{P-1} \textrm{Var}_j u_h
$$
where
$$
\textrm{Var}_j (u_h) := \sum_{n=n_h^j}^{n_h^{j+1}-1} |u_h^{n+1}-u_h^n|
$$
corresponds to the variation on $[t_h^{n_h^j},t_h^{n_h^{j+1}}[$. To study these terms, we recall that
\be{eq:unp1}
u_h^{n+1} \in \PPP_{ \frac{C(t_h^{n+1})-q_h^n}{h} } \left[u_h^n+hf_h^n\right]
\ee
by  writing the scheme in term of velocities and we state the following lemma:

\begin{lem}\label{lem:x1_moins_x0}
There exist uniformly bounded vectors $y^{n_h^j}$ ($|y^{n_h^j}|\leq \kappa_0$) such that, for all small enough $h$, for all $j\in\{0, \ldots , P\}$ and $n\in {\mathbb N} \cap [n_h^j,n_h^{j+1}[\virg $ we have
$$
x_1 \in \PPP_{ \frac{C(t_h^{n+1})-q_h^n}{h} } [x_0]\quad \Longrightarrow \quad |x_1-x_0|\leq \frac{2}{\delta}\left(|x_0-y^{n_h^j}|^2-|x_1-y^{n_h^j}|^2 \right)
$$
as soon as $x_1$ is bounded by $K$ (constant introduced in Proposition \ref{prop:uLI}). 
\end{lem}

\dem 
First for small enough $h<\nu_{min}$ and $h<\tau/2$, it comes
$$ n\in [n_h^j,n_h^{j+1}[ \Longrightarrow |t_h^{n_h^j}-t_h^n|+h<\tau/2+\tau/2 =\tau \ \textrm{ and } \ |q_h^{n_h^j}-q_h^n|\leq K|t_h^{n_h^j}-t_h^n| <r/2$$
So thanks to Lemma \ref{lem:conerentrant} (applied for $t_h^{n_h^j},q_h^{n_h^j}$), we know that there exist unit vectors $v^{n_h^j}$ such that \be{eq:vni} n\in
[n_h^j,n_h^{j+1}[\quad \Longrightarrow \quad \overline{B}(-\kappa_0 v^{n_h^j},\delta/2) \subset \frac{C(t_h^{n+1})-q_h^n}{h}. \ee Indeed $v^{n_h^j}$ is ``a good direction'', given by the admissibility assumption, associated to the point $(t_h^{n_h^j},q_h^{n_h^j})$.\\
Then, we develop similar arguments than those used in \cite{CM} and \cite{F-Aline, Paoli, DMP}. In these works, the set onto which the velocity is projected was convex. In the present case, the set $\frac{C(t_h^{n+1})-q_h^n}{h}$ is $\eta/h$-prox-regular, which is slightly weaker. \\
Let $n$ belong to $[n_h^j,n_h^{j+1}[$. We define
$$
z^{n_h^j}:=y^{n_h^j}+\frac{\delta}{2} \frac{x_0-x_1}{|x_0-x_1|}
\quad\hbox{where}\quad
y^{n_h^j}:=-\kappa_0 v^{n_h^j}.$$
(Here we suppose $x_0\neq x_1$, else the desired result is obvious.)
We have $$
z^{n_h^j}\in \overline{B}(-\kappa_0 v^{n_h^j},\delta/2)\subset
\frac{C(t_h^{n+1})-q_h^n}{h}.
$$
The point $x_1$ being the projection of $x_0$
onto the $\frac{\eta}{h}$-prox-regular set $\frac{C(t_h^{n+1})-q_h^n}{h}$, we get
$$
\langle x_0-x_1,z^{n_h^j}-x_1 \rangle \leq \frac{h}{2\eta}|x_0-x_1| |z^{n_h^j}-x_1|^2
$$
thanks to the hypomonotonicity property (see Proposition \ref{prop:hypo}).
Thus
\begin{eqnarray*}
 |x_0-y^{n_h^j}|^2&=&|x_1-y^{n_h^j}|^2+|x_0-x_1|^2+2\langle z^{n_h^j}-y^{n_h^j},x_0-x_1\rangle+2\langle x_1-z^{n_h^j},x_0-x_1\rangle\\
&\geq&|x_1-y^{n_h^j}|^2+2\langle z^{n_h^j}-y^{n_h^j},x_0-x_1\rangle - \frac{h}{\eta}|x_0-x_1| |z^{n_h^j}-x_1|^2\\
&\geq&|x_1-y^{n_h^j}|^2+\delta |x_0-x_1| - \frac{h}{\eta}|x_0-x_1| |z^{n_h^j}-x_1|^2.
\end{eqnarray*}
Using that the vectors $z^{n_h^j}$ are uniformly bounded by $\kappa_0+\delta/2$ and that $x_1$ is bounded by $K$, it follows that for $h \leq \frac{\eta \delta}{2(\kappa_0 + \delta/2 + K)^2}$
$$ |x_0-y^{n_h^j}|^2 \geq |x_1-y^{n_h^j}|^2+\frac{\delta}{2} |x_0-x_1|. $$
This, together with the fact that the vectors $y^{n_h^j}$ are uniformly bounded by $\kappa_0$, ends the proof of Lemma~\ref{lem:x1_moins_x0}. \findem

\mb We now come back to the proof of Theorem \ref{thm:uBV}. For $n$ in $[n_h^j,n_h^{j+1}[$, using~(\ref{eq:unp1}) and the previous lemma (with $x_0=u_h^n+hf_h^n$ and $x_1=u_h^{n+1}$), it comes
\begin{eqnarray*}
 |u_h^{n+1}-u_h^n-hf_h^n|&\leq&\frac{2}{\delta}\left(|x_0-y^{n_h^j}|^2-|x_1-y^{n_h^j}|^2\right)\\
&\leq&\frac{2}{\delta}\left(|u_h^n+hf_h^n-y^{n_h^j}|^2-|u_h^{n+1}-y^{n_h^j}|^2\right)\\
&\leq&\frac{2}{\delta}\left(|u_h^n-y^{n_h^j}|^2-|u_h^{n+1}-y^{n_h^j}|^2\right) +\frac{2}{\delta}|hf_h^n|^2+\frac{4}{\delta}|hf_h^n||u_h^n-y^{n_h^j}|\\
&\leq&\frac{2}{\delta}\left(|u_h^n-y^{n_h^j}|^2-|u_h^{n+1}-y^{n_h^j}|^2\right)
+\frac{2}{\delta}|hf_h^n|^2+\frac{4}{\delta}|hf_h^n|\left(K+\kappa_0\right).
\end{eqnarray*}
By summing up these terms
for $n$ from $n_h^j$ to $n_h^{j+1}-1$ we get
\begin{eqnarray*}
\textrm{Var}_j (u_h) = \sum_{n=n_h^j}^{n_h^{j+1}-1} |u_h^{n+1}-u_h^n|&\leq& \frac{2}{\delta}\left(|u_h^{n_h^{j}}-y^{n_h^j}|^2-|u_h^{n_h^{j+1}}-y^{n_h^j}|^2\right) +\sum_{n=n_h^j}^{n_h^{j+1}-1}\frac{2}{\delta}|hf_h^n|^2\\
&&\quad\quad+\frac{4}{\delta}\left(K+\kappa_0+\frac{\delta}{4}\right)\sum_{n=n_h^j}^{n_h^{j+1}-1}|hf_h^n|
\end{eqnarray*}
and finally
\begin{eqnarray*}
\textrm{Var} (u_h) =\sum_{j=0}^{P-1} \textrm{Var}_j (u_h) & \leq &
\frac{2}{\delta}\sum_{j=0}^{P-1}\left(|u_h^{n_h^{j}}-y^{n_h^j}|^2-|u_h^{n_h^{j+1}}-y^{n_h^j}|^2\right)
+ \frac{2}{\delta}\left(\int_0^T F(t) dt \right)^2\\
& &\quad\quad+\frac{4}{\delta}\left(K+\kappa_0+\frac{\delta}{4}\right)\left(\int_0^T F(t) dt \right)\\
&\leq&\frac{4}{\delta}\left(K+\kappa_0\right)^2 P+
\frac{2}{\delta} \left(\int_0^T F(t) dt \right)^2+\frac{4}{\delta}\left(K+\kappa_0+\frac{\delta}{4}\right) \left(\int_0^T F(t) dt \right).
\end{eqnarray*}
This completes the proof of
Theorem~\ref{thm:uBV}, since $P$ does not depend on $h$ from~(\ref{eq:boundP}). 
\findem

\subsection{Extraction of a convergent subsequence}  \label{subsec:32}

The previous uniform bounds on the computed velocities allow us to extract a convergent subsequence.


\begin{prop} \label{prop:q_u_cv}
There exist subsequences of $(q_h)$ and $(u_h)$  (still denoted by $q_h$ and $u_h$) which respectively converge to $q\in W^{1,\infty}(\I,\R^d)$ and $u\in BV(\I,\R^d)$. Moreover $(u_h)_h$ strongly converges to $u$ in $L^1(\I,\R^d)$ and for all $t\in \I$ 
\be{eq:CC} q(t)\in C(t). \ee
The initial condition is satisfied: $q(0)=q_0$ and $u(0)=u_0$. 
\end{prop}

\dem By Proposition \ref{prop:uLI}, the sequence $(u_h)_h$ is bounded in $L^\infty(\I,\R^d)$. 
Arzel\`a-Ascoli theorem proves that $(q_h)_h$ is also relatively compact in $W^{1,\infty}(\I,\R^d)$. So up to a subsequence, we can assume that $q_h$ strongly converges to $q\in W^{1,\infty}(\I,\R^d)$. \\
Moreover, as $(u_h)$ is bounded in $BV(\I,\R^d)$, there exists a subsequence (still denoted by $u_h$) converging to $u$ in $L^1(\I,\R^d)$. It is easy to show that necessarily $u = \dot{q}$ in the distributional sense. In addition, by the uniform bound of the variation $\textrm{Var}(u_h)$, $u$ belongs to $BV(\I,\R^d)$ with $\textrm{Var}(u) \leq \sup_{h} \textrm{Var}(u_h)$. \\
Inclusion (\ref{eq:CC}) is a direct consequence of the uniform convergence of $q_h$ to $q$ together with (\ref{eq:feasible_config2}). Let us now check the last point concerning the initial condition. Since for all $h>0$, $q_h(0)=q_0$ so $q(0)=q_0$. Moreover, the initial point $q_0\in \textrm{Int}(C(0))$. The maps $q_h$ and $q$ are Lipschitz with the same constant, which implies there exist $s>0$ and $l>0$ such that for all $t\in[0,s]$ and all small enough $h>0$
$$ d_{\partial C(t)}(q_h(t)) +d_{\partial C(t)}(q(t)) \geq l.$$
So the computed points are far away from the boundary of $C(\cdot)$ during the whole interval $[0,s]$. That implies
$$ q_h^{n+1}= 2q_h^{n}-q_h^{n-1} +h^2 f_h^{n},$$
in other words 
$$ u_h^{n+1} =u_h^n+hf_h^n$$
for $n$ verifying $t_h^{n+1}\in [0,s]$.
So we deduce that 
$$ |u_h^{n+1}-u_h^n|\leq \int_{t_h^n}^{t_h^{n+1}} F(t) dt,$$
which gives by summing all these inequalities
$$ \sup_{t\in [0,s]} \ |u_h(t)-u_0|\leq \int_{0}^{s+h} F(t) dt.$$ 
Since $u_h$ converges almost everywhere to $u$ (up to a subsequence), it comes that
$$  \|u-u_0\|_{L^\infty([0,s],\R^d)} \leq \int_{0}^{s} F(t) dt.$$ 
Taking the limit when $s$ goes to $0$ implies the desired result: $u(0)=u_0$. 
\findem


\subsection{Solution of the continuous differential inclusion} \label{subsec:33}
In this subsection, we prove that the limit function (obtained in the previous subsection) satisfies the differential inclusion of Problem (\ref{eq:includiff}), according to Definition \ref{def:solution}.


\begin{prop} The limit function $q$ satisfies the continuous differential inclusion: there exists $k\in BV(\I,\R^d)$ such that in the sense of time-measures
\be{eq2stobis} d\dot{q} + dk = f(t,q)dt  \ee
 and the differential measure $dk$ is supported on $\left\{t,\ q(t)\in\partial C(t)\right\}$~:
\be{eqK2bis}  |k|(t) = \int_0^t {\bf 1}_{q(s) \in \partial C(s)} d|k|(s), \qquad k(t)=\int_0^t \xi(s) d|k|(s),\ee
with $\xi(s)\in \NN(C(s),q(s))$, $|\xi(s)|=1$ and $|k|(t):= \textrm{Var}\, (k,[0,t])$.
\end{prop}

\mb The idea is to let $h$ goes to $0$ in the relation 
\be{eq:dis} u_h^n+hf_h^n-u_h^{n+1} \in \NN\left(\frac{C(t_h^{n+1})-q_h^n}{h},u_h^{n+1}\right). \ee
We refer the reader to \cite{LS}, \cite{Saisho} and \cite{F-Juliette} for similar reasonings in the framework of first order differential inclusions. 

\dem The scheme 
$$ u_h^{n+1}\in \PPP_{\frac{C(t_h^{n+1})-q_h^n}{h}}[u_h^n+hf_h^n]$$
implies (\ref{eq:dis}) and so
\be{aze} u_h^n+hf_h^n-u_h^{n+1} \in \NN(C(t_h^{n+1}),q_h^{n+1}).\ee
Let us define a piecewise-constant function $k_h$, defined for $t\in[t_h^{n-1},t_h^{n}[$ by
$$k_h(t):=k_h^n := u_0-u_h^n+ \sum_{i=0}^{n-1}\int_{t_h^i}^{t_h^{i+1}} f(s,q_h^i) ds.$$
So we have
\be{eq:equadiff} 
k_h^{n+1}-k_h^{n} = u_h^{n} +hf_h^n - u_h^{n+1} \in \NN(C(t_h^{n+1}),q_h^{n+1}).\ee
By Proposition \ref{prop:q_u_cv}, $u_h$ converges to $u$ in $L^1(\I,\R^d)$ and $q_h$ uniformly converges to $q$. Moreover $f$ is Lipschitz with respect to the second variable thus we deduce that $k_h$ strongly converges to some function $k\in L^1(\I,\R^d)$ verifying for almost every $t\in \I$
$$ k(t):= u_0-u(t)+ \int_0^t f(s,q(s)) ds.$$
Thanks to the uniform bounded variation of $u_h$, $k$ belongs to $BV(\I,\R^d)$ and then (\ref{eq2stobis}) holds. \\
So it remains us to check (\ref{eqK2bis}). We recall that $\Omega$ is the set of $(t,q)$ with $t\in \I$ and $q\in C(t)$. Let $\chi:\I \times \R^d \rightarrow \R^d$ be any nonnegative smooth function compactly supported in $\textrm{Int}(\Omega)$. Then
\be{eq:lim1}  0\leq \int_0^T \chi(s,q(s)) d|k|(s) \leq \liminf_{h\to 0} \int_0^T \chi(s,q_h(s)) d|k_h|(s). \ee
However, 
$$ \int_0^T \chi(s,q_h(s)) d|k_h|(s)  = \sum_{n=0}^{N-1} \int_{t_h^n}^{t_h^{n+1}} \chi(s,q_h^n+(s-t_h^n) u_h^{n+1}) (k_h^{n+1}-k_h^{n}) ds. $$
Assume that for some integer $n$, $k_h^{n+1}-k_h^{n}\neq 0$, then $q_h^n+hu_h^n+h^2f_h^n \notin C(t_h^{n+1})$. Thus, for every $t\in [t_h^n,t_h^{n+1}[$
$$ d_{\partial C(t)}(q_h^n) \leq h\left(c_0 + K+\int_0^T F(s) ds\right),$$
because $q_h^n\in C(t_h^n)$, $|u_h^n|\leq K$ (see Proposition \ref{prop:uLI}) and $c_0$ is the Lipschitz constant of $C$. Consequently, for $h$ small enough
$$ d_{\partial C(t)}(q_h^n)+h K \leq d_H(\textrm{supp}(\chi), \partial \Omega)$$
and so $\chi(s,q_h^n+(s-t_h^{n})u_h^{n+1}) =0$, for $s\in [t_h^n,t_h^{n+1}[$.
Finally, we obtain that for $h$ small enough
$$ \int_0^T \chi(s,q_h(s)) d|k_h|(s) = 0, $$
hence from (\ref{eq:lim1}), we obtain
$$ \int_0^T \chi(s,q(s)) d|k|(s)=0, $$
for any nonnegative smooth function $\chi$, compactly supported in $\textrm{Int}(\Omega)$. Taking a sequence of such functions converging (increasingly) to ${\bf 1}_{\textrm{Int}(\Omega)}$, then we have
$$ \int_0^T {\bf 1}_{\textrm{Int}(\Omega)}(s,q(s)) d|k|(s)=0, $$
or equivalently
$$ |k|_T = \int_0^T {\bf 1}_{\partial \Omega}(s,q(s)) d|k|(s) = \int_0^T {\bf 1}_{\partial C(s)}(q(s)) d|k|(s).$$
To finish, it remains us to check that ``$dk(s) \in \NN(C(s),q(s))$''. We move the reader to \cite{LS} (the end of the proof for Theorem 1.1 in \cite{LS}) for precise details. Indeed the arguments relie on the hypomonotonicity property of the proximal normal cones (see Proposition \ref{prop:hypo}). By (\ref{aze}), $u_h^n+hf_h^n-u_h^{n+1}$ belongs to $\NN(C(t_h^{n+1}),q_h^{n+1})$ for all integer $n$. So for every continuous map $\phi:\I \to \R^d$ such that $\phi(t)\in C(t)$, we have
$$ \langle \phi(t_h^{n+1})-q_h^{n+1},u_h^n+hf_h^n-u_h^{n+1}\rangle \leq \frac{1}{2\eta} |u_h^n+hf_h^n-u_h^{n+1}| |\phi(t_h^{n+1})-q_h^{n+1}|^2.$$ 
Summing all these inequalities from $0$ to $N-1$, it comes with (\ref{eq:equadiff}) for every nonnegative function $\psi$
\be{eq:disss} \int_0^T \psi(t) \langle \phi(\theta_h(t))-q_h(\theta_h(t)),dk_h(t)\rangle \leq \frac{1}{2\eta} \int_0^T  \psi(t) |\phi(\theta_h(t))-q_h(\theta_h(t))|^2 d|k_h|(t) \ee
where we denote $\theta_h(t)=t_h^{n+1}$ for $t\in[t_h^{n},t_h^{n+1}[$. Since $dk_h$ and $d|k_h|$ are uniformly bounded measures. Up to extract a subsequence, we can assume that they are weakly convergent to $dk$ (the differential measure of $k$) and $da$ (where $da$ is a nonnegative measure). Necessarily, the measure $dk$ is absolutely continuous with respect to the measure $da$. So there exists a bounded and measurable function $g$ such that $dk=g da$. Taking the limit in (\ref{eq:disss}) when $h$ goes to $0$, it comes 
$$ \int_0^T \psi(t) \langle \phi(t)-q(t),g(t)\rangle da(t) \leq \frac{1}{2\eta} \int_0^T \psi(t) |\phi(t)-q(t)|^2 da(t).$$
Since this inequality holds for every nonnegative function $\psi$, we deduce that for all $t\in \I$ and all map $\phi$
$$ \langle \phi(t)-q(t),g(t)\rangle \leq \frac{1}{2\eta} |\phi(t)-q(t)|^2,$$
which yields by Proposition \ref{prop:hypo2} that $g(t)\in \NN(C(t),q(t))$. Indeed for every $t_0\in \I$ and $\phi_0\in C(t_0)$, there exists a continuous map $\phi:\I \rightarrow \R^d$ satisfying 
$$ \left\{ \begin{array}{l}
            \phi(t)\in C(t), \quad \forall t\in\I \vsp \\
            \phi(t_0)=\phi_0. 
           \end{array} \right. $$
It suffices to consider the solution of the following sweeping process (see \cite{Thibsweep})
$$ \left\{ \begin{array}{l}
            - \dot{\phi}(t) \in \NN(C(t),\phi(t)) \quad \textrm{for a.e.}\ t\in\I \vsp \\
            \phi(t_0)=\phi_0. 
           \end{array} \right. $$
That also concludes the proof of (\ref{eqK2bis}).
\findem

\subsection{Collision law} \label{subsec:34}


\mb Finally, Theorem~\ref{thm:cv} will be proved, provided that the collision law is satisfied for the limits $u$ and $q$,
which is the aim of the following proposition.
\begin{prop}
 \label{prop:CollLaw} The impact law is satisfied:
$$ \forall t_0\in \overset{\circ}{\I} \virg u^+(t_0) = \PPP_{\C_{t_0,q(t_0)}}(u^-(t_0)).$$
\end{prop}

\dem Note that, from Proposition \ref{prop:q_u_cv}, $u\in BV(\I,\R^d)$, so that the left-sided $u^-(t_0)$ and the right-sided $u^+(t_0)$ limits are well-defined. \\
The proof is quite technical so for an easy reference, we remember the definitions of the sets $\C_{t_0,q_0}$ (see Definition \ref{def:cone})
$$ \C_{t_0,q_0}:=\left\{v =\lim_{\epsilon \to 0^+} v_\epsilon, \ \textrm{ with } \ v_\epsilon \in \frac{C(t_0+\epsilon)-q_0}{\epsilon}\right\}.$$ 
Moreover, we recall that these sets are nonempty due to Proposition \ref{prop:assumption}.
From now on, let us fix the instant $t_0\in \I$.
The desired property
\be{eq:coll} u^+(t_0) = \PPP_{\C_{t_0,q(t_0)}}(u^-(t_0)) \ee
 can be seen as the limit (for $h$ going to $0$) of the ``discretized property''
\be{eq:projete0}
u_h^{n+1} \in \PPP_{\frac{C(t_h^{n+1})-q_h^n}{h}} [u_h^n+hf^n].
 \ee

\mb { \bf Step 1:} We claim that
\be{eq:u+}
u^+(t_0)\in C_{t_0,q(t_0)}.
\ee
By definition, $q\in W^{1,\infty}(\I,\R^d)$ and $u\in BV(\I,\R^d)$ so we have that 
$$ u^+(t_0) := \lim_{\epsilon \to 0^+} u(t_0 +\epsilon) = \lim_{\epsilon \to 0^+} \frac{q(t_0+\epsilon) - q(t_0)}{\epsilon}.$$
The last equality comes from
$$ \begin{array}{rcl}
\dsp \left| u^+(t_0) -\frac{q(t_0+\alpha)-q(t_0)}{\alpha}\right| & = & \dsp \left| u^+(t_0) -\frac{1}{\alpha}\int_{t_0}^{t_0+\alpha} u(s) ds\right| \vsp \vsp \\
& \leq & \dsp \sup_{s\in[t_0,t_0+\alpha]} |u^+(t_0)-u(s)| \xrightarrow[\alpha \to 0]{} 0.
\end{array} $$
Since $q(t_0+\epsilon)\in C(t_0+\epsilon)$, 
$$ u^+(t_0) \in \lim_{\epsilon \to 0^+} \frac{C(t_0+\epsilon) - q(t_0)}{\epsilon}=\C_{t_0,q(t_0)},$$
which completes the proof of (\ref{eq:u+}).

\gb
Let us now come back to the proof of~(\ref{eq:coll}). As we just proved $u^+(t_0)\in \C_{t_0,q(t_0)}$ and since
$\C_{t_0,q(t_0)}$ is a convex set (see Proposition \ref{prop:coneconvex}), (\ref{eq:coll}) is equivalent to 
\be{amontrer} \forall w\in \C_{t_0,q(t_0)}, \qquad
\langle u^-(t_0) -u^+(t_0), w-u^+(t_0) \rangle \leq 0.\ee
So, in the following, let us fix $w\in \C_{t_0,q(t_0)}$. In Step 2, we construct a family of points $w_\nu$ for $\nu>0$ such that $w_\nu$ tends to $w$ when $\nu$ goes to zero satisfying
$w_\nu\in \frac{C(t+h)-q}{h}$ for $h$ sufficiently small and $(t,q)$ close to $(t_0,q(t_0))$. Then in Step 3, for each $\nu$, we go to the limit  on $h$, $t$ and $q$ to show that $\langle u^-(t_0)-u^+(t_0), w_\nu-u^+(t_0) \rangle \leq 0$ and finally, we make $\nu$ go to zero to conclude.

\mb{\bf Step 2}: From the admissibility assumption, there exist a neighborhood $U \subset \Omega \subset \I \times \R^d$ around $(t_0,q(t_0))$ and a ``good direction'' $\zeta \in \R^d$ such that for all $(t,q)\in U$
\begin{equation}
\forall v\in \NN(C(t),q) \virg \langle \zeta ,v \rangle \leq -\delta |v|, \label{eq:gooddir}
 \end{equation}
with a numerical constant $\delta>0$. For $\nu\in ]0,1[$, we consider the point $w_\nu:=w+\nu \zeta $. We claim that for each fixed $\nu>0$, there are $\epsilon_\nu$ and $h_\nu$ such that for every $h<h_\nu$,
$(t,q)\in U$, we have
\be{claim1} |t-t_0|+|q-q(t_0)|\leq \epsilon_\nu \Longrightarrow  w_\nu\in \frac{C(t+h)-q}{h} \ee
Let us detail this point. \\
Thanks to Proposition \ref{prop:assumption}, we know that for any $\theta$ there exists a neighborhood $V_\theta \subset U$ of $(t_0,q(t_0))$ and $h_\theta$ such that for all $(t,q)\in V_\theta$ and $h\in ]0,h_\theta[$ 
$$ d_{\frac{C(t+h)-q}{h}}(w) \leq \theta.$$
Let us denote $\tilde{w}_{t,q}$ a point of $\PPP_{\frac{C(t+h)-q}{h}}(w)$. \\
If $h\leq h_\theta \leq 2\eta/(c_0+|w|+1)$, it can be proved that for $(t,q)\in V_\theta$, $w_\nu\in \frac{C(t+h)-q}{h}$ i.e $q+hw_\nu \in C(t+h)$.
By Proposition \ref{prop:hypo2} and as $d_{C(t+h)}(q+hw_\nu)\leq 2\eta$, it suffices to show that for all $x\in C(t+h)$ and $v\in \NN(C(t+h),x)$
\be{claim1bis}
\langle q+h w_\nu -x,v\rangle \leq \frac{|v|}{2\eta}|  q+h w_\nu -x |^2 \ee
as soon as 
\be{eq:conditionsurx} |x-q+hw_\nu|\leq h(c_0+|w|+1). \ee \\
Let $x\in C(t+h)$ satisfying (\ref{eq:conditionsurx}), we have
\begin{align*}
 \langle q+h w_\nu -x,v\rangle  & = \langle  q+ hw + \nu h \zeta -x,v \rangle \\
 & =  \langle  q+ h \tilde{w}_{t,q} -x,v \rangle + \nu h \langle \zeta,v\rangle + h\langle w-\tilde{w}_{t,q} , v \rangle \\
 & \leq \frac{1}{2\eta}|v| |q+ h \tilde{w}_{t,q} -x|^2 - \nu h \delta |v| + h|w-\tilde{w}_{t,q}| |v| \\
 & \leq |v| \left[ \frac{1}{2\eta}|q+ h \tilde{w}_{t,q} -x|^2 - \nu h \delta + h \theta \right],
\end{align*}
where we have used that $q+ h \tilde{w}_{t,q}$ and $x$ belong to $C(t+h)$ with the hypomonotonicity property of the proximal normal cone and (\ref{eq:gooddir}). Then taking $\theta \leq \min\{\nu\delta/2,c_1\}$ with $c_1:=1/(6(c_0+2|w|+2))$, we get
\begin{align*} 
|q+ h \tilde{w}_{t,q} -x|^2 & \leq  |q+ h w_\nu -x|^2 + h^2|w_\nu-\tilde{w}_{t,q}|^2+2h|w_\nu-\tilde{w}_{t,q}| |q-x| \\
 & \leq |q+ h w_\nu -x|^2 + h^2
\end{align*}
and so
\begin{align*}
 \frac{1}{2\eta}|q+ h \tilde{w}_{t,q} -x|^2 - \nu h \delta + h \theta & \leq \frac{1}{2\eta}|q+ h \tilde{w}_{t,q} -x|^2 - h \frac{\nu \delta}{2} \\
& \leq \frac{1}{2\eta}|q+ h w_\nu -x|^2 - h \left[\frac{\nu \delta}{2}-\frac{h}{2\eta} \right] \\
& \leq \frac{1}{2\eta}|q+ h w_\nu -x|^2 
\end{align*}
as soon as $h\leq \nu \delta \eta$.
Thus, for $\nu\leq c_1$, $\theta \leq \min\{\nu\delta/2,c_1\}$, $h\leq h_\nu:=\min\{h_\theta,\nu\delta \eta\}$ and $(t,q)\in V_\theta$, it comes
$$  \frac{1}{2\eta}|q+ h \tilde{w}_{t,q} -x|^2 - \nu h \delta + h \theta  \leq \frac{1}{2\eta}|q+ h w_\nu -x|^2 $$ 
for all $x\in C(t+h)$ satisfying (\ref{eq:conditionsurx}) and $v\in\NN(C(t+h),x)$.
That proves (\ref{claim1bis}) and so (\ref{claim1}).

\mb {\bf Step 3}:
Let us now fix the parameter $\nu\leq c_1$. \\
Thanks to the uniform Lipschitz regularity of the maps $q_h$ and their uniform convergence towards $q$, there exists $\tilde{h}_\nu\leq h_\nu$ such that for $\epsilon\leq \epsilon_\nu/(2+2K)$ and $h\leq \tilde{h}_\nu$,
$$ t_h^k,t_h^{k+1} \in [t_0-\epsilon,t_0+\epsilon] \Longrightarrow |t_h^{k+1}-t_0|+|q_h^k-q(t_0)|\leq \epsilon_\nu.$$
We recall that 
$$ K:= \sup_h \|u_h \|_{L^\infty(\I,\R^d)}<\infty.$$
Consequently, as $q_h^{k}\in C(t_h^k)$, the property (\ref{claim1}) (with $t=t_h^k$) gives $w_\nu\in \frac{C(t_h^{k+1})-q_h^k}{h}$. Moreover (\ref{eq:projete0}) describes that 
$$ u_h^k+hf^k-u_h^{k+1}\in \NN\left(\frac{C(t_h^{k+1})-q_h^k}{h},u_h^{k+1}\right).$$
Therefore, $\frac{C(t_h^{k+1})-q_h^k}{h}$ being $\frac{\eta}{h}$-prox-regular, we have (due to the hypomonotonicity property)
\be{eq:eq}
\langle u_h^k+hf_h^k-u_h^{k+1},
w_\nu-u_h^{k+1} \rangle \leq \frac{h}{2\eta} | u_h^k+hf_h^k-u_h^{k+1}| |w_\nu-u_h^{k+1}|^2.
\ee
We sum up these inequalities for $k$ from $n$ to $p$, integers chosen such that $t_h^n$ is the first time step in $[t_0-\epsilon,t_0-\epsilon+h]$ and
$t_h^{p}$ the last one in $[t_0+\epsilon-h,t_0+\epsilon]$. First, we know that 
\be{eq:1fin} \left|\sum_{k=n}^{p} h\langle f^k,
w_\nu-u_h^{k+1} \rangle\right| \leq \left(|w|+K+1\right)\int_{t_0-\epsilon}^{t_0+\epsilon+h} F(t) dt, \ee 
with
$K:=\sup_h \|u_h \|_\infty$. We also have 
\be{eq:2fin} \sum_{k=n}^{p} \langle u_h^k-u_h^{k+1},w_\nu \rangle = \langle
u_h(t_h^{n-1})-u_h(t_h^{p}), w_\nu \rangle.\ee 
We deal with the remainder as follows: 
$$ \sum_{k=n}^{p} \langle u_h^k-u_h^{k+1}, -u_h^{k+1} \rangle = \sum_{k=n}^{p} \langle u_h^k-u_h^{k+1}, u_h^{k} \rangle - |u_h^n|^2+|u_h^{p+1}|^2,$$
which gives
\begin{align}
 \sum_{k=n}^{p} \langle u_h^k-u_h^{k+1}, -u_h^{k+1} \rangle & = \frac{1}{2} \sum_{k=n}^p |u_h^k-u_h^{k+1}|^2 + \frac{1}{2}\left[ -|u_h(t_h^{n-1})|^2+|u_h(t_h^{p})|^2\right] \nonumber \\
 & = \frac{1}{2}\textrm{Var}_2(u_h)^2_{[t_h^{n-1},t_h^p]}+ \frac{1}{2}\left[ -|u_h(t_h^{n-1})|^2+|u_h(t_h^p)|^2\right], \label{eq:3}
\end{align}
where we wrote $\textrm{Var}_2$ for the $L^2$-variation of a function. Using (\ref{eq:eq}), (\ref{eq:1fin}), (\ref{eq:2fin}) and
(\ref{eq:3}), we finally get~:
\begin{align*}
\lefteqn{ \frac{1}{2}\textrm{Var}_2(u_h)_{[t_h^{n-1},t_h^p]}^2+ \frac{1}{2}\left[ -|u_h(t_h^{n-1})|^2+|u_h(t_h^{p})|^2\right] + \langle u_h(t_h^{n-1})-u_h(t_h^{p}), w_\nu \rangle} & & \\
& &  \leq (|w|+K+1) \int_{t_0-\epsilon}^{t_0+\epsilon+h} F(t) dt + \frac{h}{2\eta} \sum_{k=n}^p | u_h^k+hf_h^k-u_h^{k+1}| |w_\nu-u_h^{k+1}|^2.
\end{align*}
However
$$ \sum_{k=n}^p | u_h^k+hf_h^k-u_h^{k+1}| \leq \sum_{k=n}^p | u_h^k-u_h^{k+1}| + \sum_{k=n}^p h|f_h^k| \leq \textrm{Var}(u_h) + \int_0^T F(t)dt \leq B_1$$
for some numerical constant, due to Theorem \ref{thm:uBV} and
$$ |w_\nu-u_h^{k+1}|^2 \leq (|w| + 1+ K)^2 =B_2^2,$$
due to Proposition \ref{prop:uLI}.
Consequently, we deduce that
\begin{align}
 & \frac{1}{2}\textrm{Var}_2(u_h)_{[t_h^{n-1},t_h^p]}^2+ \frac{1}{2}\left[ -|u_h(t_h^{n-1})|^2+|u_h(t_h^{p})|^2\right] + \langle u_h(t_h^{n-1})-u_h(t_h^{p}), w_\nu \rangle \nonumber \\
&  \hspace{5cm} \leq B_2 \int_{t_0-\epsilon}^{t_0+\epsilon+h} F(t) dt + \frac{h}{2\eta} B_1 B_2^2. \label{eq:b1b2}
\end{align}
Let us now choose a sequence of $\epsilon_m$ going to zero, such that $u_h$ pointwisely converges to $u$ at the instants $t_0-\epsilon_m$ and $t_0+\epsilon_m$ (which is possible as $u_h$ converges almost everywhere towards $u$). For each $\epsilon_m$ and $h\leq \tilde{h}_\nu$, we have shown that inequality (\ref{eq:b1b2}) holds. Then, passing to the limit for $h\to 0$ we get
\begin{align*} \frac{1}{2}\textrm{Var}_2(u)_{[t_0-\epsilon_m,t_0+\epsilon_m]}^2+ \frac{1}{2}\left[ -|u(t_0-\epsilon_m)|^2+|u(t_0+\epsilon_m)|^2\right] & \\
 & \hspace{-4cm}+\langle u(t_0-\epsilon_m)-u(t_0+\epsilon_m), w_\nu \rangle\leq B_2 \int_{t_0-\epsilon_m}^{t_0+\epsilon_m} F(t) dt,
\end{align*}
which gives for $\epsilon_m\to 0$
$$ \frac{1}{2}\textrm{Var}_2(u)_{[t_0^-,t_0^+]}^2+ \frac{1}{2}\left[ -|u^-(t_0)|^2+|u^+(t_0)|^2\right] + \langle u^-(t_0)-u^+(t_0), w_\nu \rangle\leq 0.$$
Finally we obtain
$$ \frac{1}{2}\left| u^{+}(t_0) - u^-(t_0)\right|^2 + \frac{1}{2}\left[ -|u^-(t_0)|^2+|u^+(t_0)|^2\right] + \langle u^-(t_0)-u^+(t_0), w_\nu \rangle\leq 0.$$
By expanding the square quantities, we obtain for all $\nu<c_1$
\be{eq:amontrer_nu} \langle u^-(t_0)-u^+(t_0), w_\nu-u^+(t_0) \rangle \leq 0.\ee
Recall that $w_\nu=w+\nu \zeta$, we obtain (\ref{amontrer}) by letting $\nu$ go to $0$ in~(\ref{eq:amontrer_nu}). \findem

\section{A particular case} \label{sec:kpart}

As explained in the introduction, all the second-order differential inclusions, already studied in the literature, concern a particular case where the moving set $C(\cdot)$ is given by a finite number of constraints. This section is also devoted to prove that the previous abstract result covers this case, as soon as the constraints satisfy some reasonable assumptions. 

\mb So we consider the Euclidean space $\R^d$, $B$ the closed unit ball in $\R^{d+1}$ and $\I=[0,T]$ a bounded time-interval. For $i\in\{1,...,p\}$ let $g_i:\ \I \times \R^d \rightarrow \R$ be functions (which can be thought as ``constraints''). For $t\in\I$, we introduce the sets
$$ Q_i(t):=\left\{x\in \R^d,\ g_i(t,x)\geq 0\right\},$$
and the following one
$$ Q(t):= \bigcap_{i=1}^p Q_i(t),$$
which represents the set of ``feasible configurations $x$''.
We remember that $\Omega:=\{(t,x),\ t\in\I, \ x\in Q(t)\}$ and we similarly define $\Omega_i$ in replacing $Q(t)$ by $Q_i(t)$. Moreover we assume that there exist $\alpha,\beta,M,\kappa>0$ such that $g_i\in C^2\left(\Omega +\kappa B \right)$ and satisfies~:
\begin{equation} \label{gradg}
\forall (t,x)\in \Omega_i+\kappa B, \qquad \alpha \leq |\nabla_x g_i(t,x)| \leq \beta, \tag{$A1$}
\end{equation}
\begin{equation} \label{dtg}
\forall (t,x)\in \Omega_i+\kappa B, \qquad |\partial_t g_i(t,x)| \leq \beta, \tag{$A2$}
\end{equation}
\begin{equation} \label{hessg}
\forall (t,x)\in \Omega_i+\kappa B, \qquad |D_{x}^2 g_i(t,x)| \leq M \tag{$A3$}
\end{equation}
and
\begin{equation} \label{hesstg}
\forall (t,x)\in \Omega_i+\kappa B, \qquad |\partial_t^2 g_i(t,x)| + |\partial_t \nabla_x g_i(t,x)| \leq M. \tag{$A4$}
\end{equation}
For all $t\in\I$, we denote by
$$ I(t,x):=\left\{i,\ g_i(t,x)=0\right\} $$ the set of ``active contraints'' and for $\rho>0$ 
$$ I_\rho(t,x):=\left\{i,\ g_i(t,x)\leq \rho\right\}.$$
We suppose that there exist constants $\rho,\gamma>0$ such that for all $x\in Q(t)$ and all nonnegative reals $\lambda_i$
\be{reverse}
\sum_{i\in I_\rho(t,x)} \lambda_i |\nabla g_i(t,x)| \leq \gamma \left| \sum_{i\in I_\rho(t,x)} \lambda_i \nabla g_i(t,x)\right |, \tag{$R_\rho$}
\ee

\begin{thm} \label{thm:proxc} Under the assumptions (\ref{gradg}), (\ref{hessg}) and ($R_0$), there exists $\eta:=\eta(\alpha,M,\gamma)$ such that the set $Q(t)$ is $\eta$-prox-regular for all $t\in\I$. 
\end{thm}

\dem The time $t$ is fixed in this proof so for simplicity we omit it in the notations.
We will follow the arguments and the ideas of \cite{TheseJu} (Subsections 3.1 and 3.2) and \cite{Juliette-M} (Subsection 2.2) where the desired result is already proved in the case of convex constraints $g_i$. \\
First let us study the set $Q_i$ for a fixed index $i$. We refer the reader to Proposition 3.2 of \cite{TheseJu} for the following well-known fact. Due to the assumptions, $Q_i$ has a $C^1$-boundary 
$$ \partial Q_i = \left\{x\in \R^d, \ g_i(x)=0\right\}$$
and for $x\in \partial Q_i$, its proximal normal cone is given by
$$ \NN(Q_i,x) = -\R^+ \nabla g_i(x).$$
We now want to check that $Q_i$ is uniformly prox-regular with a constant $\eta_0$. It suffices to check the hypomonotonicity property: for all $x\in \partial Q_i$ and $y\in Q_i$
\be{eq:hypo} \langle x-y, -\nabla g_i(x) \rangle \geq -\frac{1}{2\eta_0}|x-y|^2 |\nabla g_i(x)|. \ee
Indeed since $g_i(x)=0$, a first order expansion together with Assumption (\ref{hessg}) give
$$ 0\leq g_i(y) \leq   \langle x-y, -\nabla g_i(x) \rangle + \frac{M}{2} |x-y|^2.$$
So (\ref{eq:hypo}) is satisfied with $\eta_0:=\alpha/M$, hence $Q_i$ is $\eta_0$-prox-regular (thanks to Proposition \ref{prop:hypo}). \\
Then let us study $Q$ the intersection of sets $Q_i$. We first have to prove that for all $x\in \partial Q$
\be{eq:cone} \NN(Q,x) = \sum_{i\in I(x)} \NN(Q_i,x) = -\sum_{i\in I(x)} \R^+ \nabla g_i(x). \ee
Let us denote 
$${\mathcal N}_x:=-\sum_{i\in I(x)} \R^+ \nabla g_i(x).$$
The inclusion ${\mathcal N}_x \subset \NN(Q,x)$ is proved in Proposition 2.16 \cite{Juliette-M} (this part did not use the convexity of the functions $g_i$) and so we just deal with the other one.
By the way, we point out that ${\mathcal N}_x$ is the polar cone of 
$$ \Upsilon_x:=\left\{z\in \R^d, \ \forall i \in I(x), \ \langle \nabla g_i(x),z\rangle \geq 0 \right\}.$$
So using the orthogonal decomposition related to polar cones (see \cite{Moreau}), any $v \in \NN(Q,x) $ can be written $v=w +z=\PPP_{{\mathcal N}_x} v+\PPP_{\Upsilon_x}v$, with  $w \bot z$.
Suppose  $z \neq 0$. Since $v \in \NN(Q,x) $, there exists $t>0$ such that $ x \in \PPP_{Q}(x +tv)$. Let
$$s=\min(t,\tau) \textmd{   with   }\dsp  \tau := \min_{i
  \notin I(x)} \frac{g_i(x)}{(2\beta+\delta \alpha) |z|},$$ 
 by the well-known property of the projection, the inequality $s \leq t$ implies 
 \be{eq:pro} x \in
  \PPP_{Q}(x +s v). \ee 
From Lemma 5.2 in \cite{F-Juliette}, we know that $Q$ satisfies the second property of the admissibility: there exist a bounded covering of $\partial Q$ with balls $B((x_p,r))_p$, a collection of ``good direction $(u_p)_p$'' and constants $\rho,\delta$ such that for all $x\in B(x_p,2r)$ and all $v\in \NN(Q,x)$ 
$$ \langle v,u_p \rangle \geq \delta |v|.$$
There exists $p$ such that $x\in B(x_p,r)$ and we set
$$\tilde{x}:=x + sv-sw-\epsilon s|z|u_p=x+sz-\epsilon s|z|u_p,$$
where $\epsilon$ will be later chosen small enough ($\epsilon<<1$). We claim that $ \tilde{x}
  \in Q$. Indeed thanks to Assumption (\ref{hessg}) we have for all $i$, 
$$ g_i(\tilde{x})\geq g_i(x)+s
  \langle \nabla g_i(x), z-\epsilon |z|u_p \rangle - s^2M(1+\epsilon)^2|z|^2. $$  Consequently,
  $$\forall i\in I(x), \    g_i(\tilde{x})\geq  s\epsilon |z|\delta\alpha - s^2M(1+\epsilon)^2|z|^2  \geq 0$$ 
if $s\leq \frac{\epsilon \delta \alpha}{|z|M(1+\epsilon)^2}$. \\
Furthermore, if $i \notin  I(x) $, then $ \dsp s \leq \tau \leq
  \frac{g_i(x)}{(\epsilon\delta \alpha +(1+\epsilon) \beta) |z|}.$ Hence $$
  g_i(\tilde{x})\geq g_i(x)-s\beta|z| (1+\epsilon) \geq 0   .$$
That is why $\tilde{x} \in Q$ as soon as $s\leq \frac{\epsilon \delta \alpha}{|z|M(1+\epsilon)^2}$ and in this case $d_{Q}(x +sv) \leq |x +s
  v - \tilde{x} | = s\Big| w+\epsilon |z|u_p \Big| :=s\sqrt{A}$. By expanding
 \begin{align*}
 A & = \Big| w+\epsilon |z|u_p \Big| ^2 \\
  & = |w|^2 + \epsilon ^2 |z|^2 + 2 \epsilon |z|\langle w, u_p \rangle \\
  & \leq |w|^2 + \epsilon ^2 |z|^2 + 2 \epsilon |z||w| \\
  & \leq |w|^2 + \frac{1}{2} |z|^2 \\
  & < |v|^2 := |w|^2 + |z|^2,
  \end{align*}
  for $\epsilon\leq \frac{1}{2}\frac{|z|}{4|w|+|z|}$. Finally we have obtained that for small enough $\epsilon$ and $s$
  $$ d_{Q}(x +sv)<s|v|$$
  which leads to a contradiction with (\ref{eq:pro}). So we conclude that $z=0$, which completes the proof of (\ref{eq:cone}).\\
Finally, the prox-regularity of the set $Q$ is shown by invoking the ``reverse triangle inequality'' Assumption (\ref{reverse}), as done in Proposition 2.17 \cite{Juliette-M}.
\findem

\mb 
\begin{rem} \label{rem:conetn} Note that for all $t\in\I$, $x\in Q(t)$
$$ \Upsilon_{t,x}:=\left\{z\in \R^d, \ \forall i \in I(t,x), \ \langle \nabla_x g_i(t,x),z\rangle \geq 0 \right\} = \TT_{Q(t)}(x)$$
Indeed $\Upsilon_{t,x}$ and $\TT_{Q(t)}(x)$ are two convex cones whose polar cones are equal:
$$ {\mathcal N}_{t,x}:=-\sum_{i\in I(t,x)} \R^+ \nabla_x g_i(t,x) = \NN(Q(t),x).$$
\end{rem}

\begin{prop} Let $Q$ be the set-valued map defined at the beginning of the current section. The set $\Omega:=\left\{(t,x)\in \I \times \R^{d},\ x\in Q(t)\right\}$ is uniformly prox-regular.
Moreover for all $(t,x)\in\Omega$ with $t\in \overset{\circ}{\I}$, the set ${\mathcal C}_{t,x}$ (defined in Definition \ref{def:cone}) verifies
$$ {\mathcal C}_{t,x}= \left\{ z\in\R^d,\ \forall i\in I(t,x),\ \partial_t g_i(t,x)+\langle z,\nabla_x g_i(t,x) \rangle \geq 0 \right\}.$$
\end{prop}

\dem The set $\Omega$ is given by the functions $h_i$ as 
$$ \Omega:= \left\{(t,x)\in\R^{d+1}, \ \forall i\in\{0,...,p+1\},\  h_i(t,x)\geq 0,\ \right\}$$
with $h_i:=g_i$ for $i\in\{1,...,p\}$, $h_0(t,x):=t$ and $h_{p+1}(t,x):=T-t$.
In order to apply Theorem \ref{thm:proxc}, we check the different assumptions relatively to $h_i$: for all $i\in\{0,...,p+1\}$ 
\begin{equation} \label{nablah}
\forall (t,x)\in \Omega_i+\kappa B, \qquad \min\{\alpha,1\} \leq |\nabla_{(t,x)} h_i(t,x)| \leq \max\{2\beta,1\}, 
\end{equation}
\begin{equation} \label{hessh}
\forall (t,x)\in \Omega_i+\kappa B, \qquad |D_{(t,x)}^2 h_i(t,x)| \leq 2M,
\end{equation}
and
\be{reverseh}
\sum_{i\in I(t,x)} \lambda_i |\nabla_{(t,x)} h_i(t,x)| \leq \gamma\left(1+\frac{\beta+2}{\alpha}\right) \left| \sum_{i\in I(t,x)} \lambda_i \nabla_{(t,x)} h_i(t,x)\right |.
\ee
Indeed Properties (\ref{nablah}) and (\ref{hessh}) are obvious for $i=0$ and $i=p+1$. Let $i\in\{1,...,p\}$ and $(t,x)\in \Omega_i+\kappa B$.
Inequality (\ref{nablah}) is proved by
$$ \alpha \leq |\nabla_x g_i(t,x)| \leq |\nabla_{(t,x)} h_i(t,x)| \leq |\nabla_x g_i(t,x)|+|\partial_t g_i(t,x)| \leq 2\beta.$$
The second one (\ref{hessh}) is due to Assumptions (\ref{hessg}) and (\ref{hesstg}):
$$ |D_{(t,x)}^2 h_i(t,x)| \leq \left(|\partial_t^2 g_i(t,x)|^2 + |D_{x}^2 g_i(t,x)|^2 + 2|\partial_t \nabla_x h_i(t,x)|^2\right)^{1/2} \leq 2M.$$
Concerning Assumption (\ref{reverseh}), we just deal with the case where $I(t,x) \cap \{1,...,p\} \neq \emptyset$ (else the inequality is obvious). For all $i\in\{1,...,p\}$, we have $|\partial_t\, g_i(t,x)|\leq \frac{\beta}{\alpha} |\nabla_x\, g_i(t,x)|$ which implies
$$ |\nabla_{(t,x)}\, h_i(t,x) | \leq \left(|\partial_t\, g_i(t,x)|^2 + |\nabla_x\, g_i(t,x)|^2\right)^{1/2} \leq \left(1+\frac{\beta}{\alpha}\right) |\nabla_x\, g_i(t,x)|.$$
For $i=0$ or $i=p+1$, it comes similarly
$$ |\nabla_{(t,x)} h_i(t,x) | =1 \leq \frac{1}{\alpha} |\nabla_x\, g_{i_0}(t,x)|$$
where $i_0\in I(t,x) \cap \{1,...,p\}$.
Then (\ref{reverseh}) is involved by Assumption (\ref{reverse}).
By Theorem \ref{thm:proxc}, we conclude to the uniform prox-regularity of $\Omega$. \\
Thanks to Proposition \ref{prop:conebis}, we know that for all $(t,x)\in \Omega$ with $t\in \overset{\circ}{\I}$
$$ {\mathcal C}_{t,x} = \left\{u\in\R^d,\ (1,u) \in \TT^D_{\Omega}((t,x)) \right\} = \left\{u\in\R^d,\ (1,u) \in \TT_{\Omega}((t,x)) \right\}$$
because $\Omega$ is uniformly prox-regular (see Corollary 6.30 \cite{RW}).
By Remark \ref{rem:conetn}, we get
\begin{align*}
 \TT_{\Omega}((t,x)) & = \left\{z\in \R^{d+1}, \ \forall i \in I(t,x), \ \langle \nabla_{(t,x)}\, h_i(t,x),z\rangle \geq 0 \right\} \\
 & = \left\{z=(z_1,z_u)\in \R \times \R^{d}, \ \forall i \in I(t,x)\cap \{1,...,p\}, \ z_1 \partial_t g_i(t,x) + \langle \nabla_{x}\, g_i(t,x),z_u \rangle \geq 0 \right\},
\end{align*}
since $t\in \overset{\circ}{\I}$. Consequently, $ u\in {\mathcal C}_{t,x}$ if and only if $\partial_t g_i(t,x)+\langle u,\nabla_x\, g_i(t,x) \rangle \geq 0$ for all $i\in I(t,x)=  I(t,x)\cap \{1,...,p\} $.
\findem

\begin{prop} \label{prop:fin} Under the previous assumptions, the set-valued map $Q(\cdot)$ is admissible. More precisely, for all $(t_0,q_0)\in\Omega$ there exist a ``good direction $u$'' and constants $r,\delta$ such that for all $(t,q)\in\Omega \cap B((t_0,q_0),2r)$ and all proximal vector $v\in \NN(Q(t),q)=-\sum_{i\in I(t,q)} \R^{+} \nabla g_i(t,q) $ 
\be{es:goodd} \langle v,u \rangle \geq \delta |v|. \ee
Moreover the set-valued map $Q$ is Lipschitz continuous on $\I$.
 \end{prop}

\mb We refer the reader to Lemma 5.2 in \cite{F-Juliette} for the admissibility property and to Proposition 2.11 of \cite{Juliette} for the Lipschitz continuity.

\mb The case of a set-valued map $Q$ has already been studied in \cite{F-Aline}. We look for explaining that the current results of existence solutions for differential inclusions covers the one obtained in \cite{F-Aline}.

\begin{prop} Under the above assumptions, Problem (\ref{eq:includiff}) is equivalent to the following one~: find $q\in W^{1,\infty}(\I,\R^d)\virg \dot q\in BV(\I,\R^d)$ and time-measures $\lambda_i \in {\mathcal
  M}_+(\I)$ such that
\be{Pb:cont2}
\left\{
\begin{array}{l}
\dsp \forall t\in \I,\quad q(t)\in Q(t) \vsp \\
\dsp d\dot{q} = f(t,q)dt + \sum_{i=1}^p  \nabla_q \, g_i (t,q) d\lambda_{i}  \vsp \\
\dsp \hbox{supp}(\lambda_{i})\subset \{t\virg g_{i}(t,q(t))=0\}\hbox{ for all } i \vsp\\
\dsp \forall t\in \I,\quad \dot q(t^+) = \PPP_{C_{t,q(t)}}\dot q(t^-)\vsp \\
\dsp q(0)=q_0 \vsp \\
\dot q (0)=u_0.
\end{array}
\right.
\ee
\end{prop}

\dem By the characterization of the proximal normal cones (see Remark \ref{rem:conetn}), it is obvious that a solution of Problem (\ref{Pb:cont2}) is a solution of Problem (\ref{eq:includiff}) too. Indeed, consider $q,\lambda_1,...,\lambda_p$ a solution of Problem (\ref{Pb:cont2}). Then we define the measure $\lambda=\lambda_1+...+\lambda_p$ in order that each measure $\lambda_i$ is absolutely continuous with respect to $\lambda$. There also exist bounded nonnegative and mesurable functions $\ell_i$ such that for all $i=1,...,p$: $d\lambda_i=\ell_i d\lambda$. Then
$$ d\dot{q} = f(t,q)dt + \sum_{i=1}^p  \nabla_q \, g_i (t,q) d\lambda_{i} = f(t,q) dt + \left(\sum_{i=1}^p  \nabla_q \, g_i (t,q) \ell_i(t) \right) d\lambda.$$
By writing $\tilde{\xi}(t) = \sum_{i=1}^p  \nabla_q \, g_i (t,q) \ell_i(t) \in -\NN(Q(t),q(t))$, we set $\xi(t)=\tilde{\xi}(t)/|\tilde{\xi}(t)|$ if $\tilde{\xi}(t)$ is non-vanishing and $\xi(t)=0$ else. So, we obtain
$$ d\dot{q} = f(t,q)dt +\xi(t) d\nu(t)$$
with $d\nu = |\tilde{\xi}(t)| d\lambda$. Since the functions $\ell_i$ are bounded by $1$, the function $\tilde{\xi}$ is bounded on $\I$ and $\nu \in {\mathcal M}_+(\I)$. We conclude that $q$ is a solution of Problem (\ref{eq:includiff}). \\
Let us explain the other relation. Let $q$ be a solution of (\ref{eq:includiff}). By definition, there exist a function $k\in BV(\I,\R^d)$ and a measurable map $\xi:\I \rightarrow \R^d$ such that 
$$ d\dot{q} + dk = f(t,q)dt $$
and
$$ k(t)=\int_0^t \xi(s) d|k|(s)$$
with $\xi(s)\in \NN(C(s),q(s))$ and $|\xi(s)|=1$.
The technical difficulty is to represent the map $\xi$ in terms of Kuhn-Tucker multipliers (this corresponds to a problem of selection for a particular set-valued map). For all $t\in\I$, we define
$$ {\mathcal F}(t):=\left\{ (\lambda_1,...,\lambda_p)\in (\R^+)^p,\ \xi(t) = -\sum_{i=1}^p \lambda_i \nabla_q g_i(t,q(t)),\ \lambda_i\neq 0 \textrm{ only if } g_i(t,q(t))=0\right\}.$$
We are now looking to obtain a measurable selection of the set-valued map ${\mathcal F}$. It is easy to see that ${\mathcal F}$ takes non-empty closed values. Moreover its graph $\Gamma_{\mathcal F}$ 
$$  \Gamma_{\mathcal F}:= \left\{ (t,\lambda_1,...,\lambda_p) \in \I \times (\R^+)^p, \ (\lambda_1,...,\lambda_p) \in {\mathcal F}(t) \right\}$$
is given by
$$ \Gamma_{\mathcal F}:= {\mathcal G}^{-1}(\{0\}) \cap \left(\I \times (\R^+)^p\right)$$
where ${\mathcal G}$ is defined as follows
$$ {\mathcal G}(t,\lambda_1,...,\lambda_p) := \xi(t)+\sum_{i\in I(t,q(t))} \lambda_i \nabla_q g_i(t,q(t)) + \left(\sum_{i\notin I(t,q(t))} |\lambda_i| \chi_{I(t,q(t))^c}(i)\right) e$$
where $e$ is any vector $e\in\R^d\setminus\{0\}$ and $\chi_{I(t,q(t))^c}$ is the characteristic function:
$$ \chi_{I(t,q(t)^c}(i)=\left\{ \begin{array}{ll}
                               0, &  i\in I(t,q(t)) \\
                               + \infty, &  i\notin I(t,q(t)). 
                              \end{array} \right.
$$
Consequently, since $\xi$ is measurable and $t\rightarrow I(t,q(t))$ is upper semicontinuous, we deduce that ${\mathcal G}$ is measurable and so it follows that the graph $\Gamma_{\mathcal F}$ is a measurable set. Then, Theorem 8.1.4 in \cite{AF} yields that the set-valued map ${\mathcal F}$ is measurable and that ${\mathcal F}$ admits a measurable selection. Let us write $(a_1,...,a_p)$ for such a measurable selection of ${\mathcal F}$. Then, we let the reader to check that the measures $d\lambda_i := a_i d|k|$ are solutions of Problem (\ref{Pb:cont2}) (since the reverse triangle inequality (\ref{reverse}), it comes that the functions $a_i$ belong to $L^\infty(\I,\R)$).
\findem

\mb We have checked that the set-valued map $Q$ takes uniformly prox-regular values (see Theorem \ref{thm:proxc}). Moreover the Lipschitz regularity and the admissibility property have already been proved (see Proposition \ref{prop:fin}). We can now apply Theorem \ref{thm:cv} and we get the following one:

\begin{thm}  \label{thm:particular} Under the above assumptions ((\ref{gradg})-(\ref{reverse})) and with $f$ satisfying (\ref{flip}) and (\ref{lingro}), the problems (\ref{Pb:2}) and (\ref{Pb:cont2}) are equivalent and admit at least one solution on any time interval $\I$.
\end{thm}

\mb We would like to finish this section with a local existence result for a constant set $Q(t)=Q$ defined by time-independent constraints.

\begin{thm} \label{thm:partlocal} Let $\kappa>0$ and consider $ Q:= \bigcap_{i=1}^p \left\{x\in \R^d,\ g_i(x)\geq 0\right\}$, where
$g_i\in C^2(Q+\kappa \overline{B}(0,1))$  satisfy (\ref{gradg}), (\ref{hessg}) and ($R_0$). Then for all initial data $q_0\in\textrm{Int}(Q)$ and $u_0\in\R^d$, there exist $T_0:=T_0(|u_0|)>0$ and $q$ solution of (\ref{Pb:}) on $[0,T_0]$ where $f$ satisfies (\ref{flip}) and (\ref{lingro}) with $ F \in L^\infty(\I)$.
\end{thm}

\dem First, $Q$ is uniformly prox-regular by Theorem \ref{thm:proxc}. Thanks to Theorem \ref{thm:appendix}, there exists $T_0:=T_0(|u_0|)>0$ such that the computed velocities $(u_h)_h$ are uniformly bounded in $L^\infty([0,T_0],\R^d)$ (for $h$ going to $0$). We can also conclude that all the computed solutions $q_h$ are uniformly bounded: there exists $L>0$ such that for all small enough $h>0$ and all $t\in[0,T_0]$~:
$$ |q_h(t)|\leq L.$$
Due to Lemma \ref{lem:r0rrho}, there exist $\rho:=\rho(L)$ and $\gamma:=\gamma(L)>0$ such that (\ref{reverse}) holds for every $q\in \overline{B}(0,2L)$. By Lemma 5.2 in \cite{F-Juliette}, the set $Q$ verifies a ``local admissibility'' property: there exist $\delta,r>0$, sequences $(x_p)_p$ and $(u_p)_p$ with $|u_p|=1$ and $x_p\in Q \cap \overline{B}(0,2L)$ such that $(B(x_p,r))_p$ is a bounded covering of the boundary $\partial Q \cap \overline{B}(0,2L)$ and 
$$ \forall p,\ \forall x\in \partial Q \cap \overline{B}(0,2L) \cap B(x_p,2r),\ \forall v\in \NN(Q,x), \quad \langle v,u_p\rangle \geq \delta |v|.$$
We let the reader to check that the proof of Theorem \ref{thm:cv} still holds since we consider points $q_h^n \in Q \cap \overline{B}(0,L)$ (indeed we can choose $r\leq L$ in order that for all $q_h^n\in \overline{B}(0,L)$, the whole ball $\overline{B}(q_h^n,r)$ is included in $B(0,2L)$ where the admissibility property is verified). Thus we get the existence of a solution on $[0,T_0]$.
\findem

\begin{lem} \label{lem:r0rrho} Let $L>0$ then there exist $\rho:=\rho(L)$ and $\gamma:=\gamma(L)>0$ such that (\ref{reverse}) holds for every $q\in \overline{B}(0,2L)$.
\end{lem}

\dem Assume that $(R_\rho)$ does not hold for all $\rho>0$ in $\overline{B}(0,2L)$. So there exist a sequence $(q_n)_n\in \overline{B}(0,2L)$ and $(\lambda_i^n)_{i,n}$ such that
$$ \sum_{i\in I_{1/n}(q_n)} \lambda_i^n |\nabla g_i(q_n)| \geq n \left|\sum_{i\in I_{1/n}(q_n)} \lambda_i^n \nabla g_i(q_n)\right|,$$
with $\lambda_i^n\geq 0$. Without loss of generality, we can assume that $\sum_{i\in I_{1/n}(q_n)} \lambda_i^n =1$. \\
By compactness, up to a subsequence, there exist $q\in \overline{B}(0,2L)$ and $(\lambda_i)_{i=1,...,p}$ such that $q_n\rightarrow q$ and for all $i\in\{1,...,p\}$, $\lambda_i^n \rightarrow \lambda_i$.
It comes with (\ref{hessg})
$$
 n\left|\sum_{i\in I_{1/n}(q_n)} \lambda_i^n \nabla g_i(q) \right| - n M|q_n-q| \leq \sum_{i\in I_{1/n}(q_n)} \lambda_i^n |\nabla g_i(q)| + M|q_n-q|.$$
So
\begin{align*}
 \left|\sum_{i\in I_{1/n}(q_n)} \lambda_i^n \nabla g_i(q) \right| & \leq \frac{1}{n}\sum_{i\in I_{1/n}(q_n)} \lambda_i^n |\nabla g_i(q)| + \left( M + \frac{M}{n}\right)|q_n-q| \\
 & \leq \frac{\beta}{n} + \left( M + \frac{M}{n}\right)|q_n-q|.
\end{align*}
Furthermore, since $g_i\in C^1(\R^d)$, $g_i$ is uniformly continuous and so it is easy to check that $I_{1/n}(q_n) \subset I(q)$ for $n$ large enough. For such $n$, ($R_0$) and (\ref{gradg}) imply
$$ \frac{\alpha}{\gamma} \leq \frac{1}{\gamma} \sum_{i\in I_{1/n}(q_n)} \lambda_i^n |\nabla g_i(q)| \leq \left|\sum_{i\in I_{1/n}(q_n)} \lambda_i^n \nabla g_i(q)\right|.$$
Finally, we obtain
$$ \frac{\alpha}{\gamma} \leq \frac{\beta}{n} + \left( M + \frac{M}{n}\right)|q_n-q|$$
which leads to a contradiction for $n\to \infty$.
\findem

\mb In Theorem \ref{thm:partlocal}, the assumptions on the gradients $(\nabla g_i)_i$ are weakened with respect to the existing results (see the introduction). Indeed the gradients of active constraints are usually supposed to be linearly independent, which implies that Assumption ($R_0$) holds on any compact (see Remark \ref{rem:indR0}). Assumption ($R_0$) permits also to deal with a large number of active constraints, while the linearly independence requires that $|I(q)|\leq d$. In addition, we emphasize that the classical geometrical assumption: for all $i\neq i'\in I(q)$ 
$$ \langle \nabla g_i(q), \nabla g_{i'}(q) \rangle \leq 0 $$
is not required.

\begin{rem} \label{rem:indR0} Let $K\subset \R^d$ be a compact set and assume that for all $q\in Q\cap K$, $(\nabla g_i(q))_{i\in I(q)}$ is linearly independent. This collection can be extended to a basis of $\R^d$, denoted by $(e_i^q)_{i=1,...,d}$. Let us consider $(\ell_j^q)_{j=1,...,d}$ the associated dual basis such that for all $i,j\in\{1,..,d\}$
$$ \ell_j^q (e_i^q) = \left\{\begin{array}{ll}
                        1, &  i=j \\
                        0, & i\neq j .
                        \end{array} \right. $$
Then for $q\in Q\cap K$ and nonnegative reals $(\lambda_i)_{i\in I(q)}$, it comes
\begin{align*} 
\sum_{i\in I(q)} \lambda_i |\nabla g_i(q)| & \leq \beta \sum_{i\in I(q)} \lambda_i \ell_i^q(\nabla g_i(q)) \leq \beta \sum_{i\in I(q)} \sum_{j=1}^d  \lambda_i \ell_j^q (\nabla g_i(q)) \\
 & \leq \beta \sum_{j=1}^d \ell_j^q\left(\sum_{i\in I(q)}\lambda_i \nabla g_i(q) \right) 
 \leq \beta \left\| \sum_{j=1}^d \ell_j^q \right\| \left|\sum_{i\in I(q)} \lambda_i \nabla g_i(q) \right|.
\end{align*}
So, we obtain ($R_0$) on $K$ with
$$ \gamma:=  \beta \sup_{q\in Q\cap K}  \left\| \sum_{j=1}^d \ell_j^q \right\| <\infty$$ by Lemma 7 of \cite{Paoli2}.
\end{rem}

\section{Appendix} \label{sec:local}

This section is devoted to another proof of uniform bounds for the computed velocities in $L^\infty(\I,\R^d)$ (Proposition \ref{prop:uLI}). However, we are looking for a proof in the framework of a constant set $C$, requiring only the prox-regularity property on $C$ and we do not assume it is admissible. We will see that we can just obtain a local result: the velocities are bounded on a time interval $[0,T_0]$, where the time $T_0$ depends on the initial conditions. It is not clear how can we extend the proof to a global result without extra properties such the admissibility.

\begin{thm} \label{thm:appendix} Assume that $C(\cdot):=C$ where $C$ is a uniformly prox-regular (possibly not admissible) set and $f$ satisfies (\ref{flip}) and (\ref{lingro}) with $ F \in L^\infty(\I)$. Then, the computed velocities (still by using the scheme (\ref{eq:qn+1})) $(u_h)_{h\leq 1}$ are uniformly bounded in $L^\infty(\I,\R^d)$ on $[0,T_{0}]$, where $T_0$ depends on the initial data $u_0$ as follows:
$$ T_0= \frac{1}{2(J+1) \left(2|u_0|+3\|F\|_\infty + \sqrt{\|F\|_\infty}\right) },$$
with a numerical constant $J$ (defined in Lemma \ref{lem:trinome}).
\end{thm}

\mb We first detail three technical lemmas and postpone the proof after them.

\begin{lem} \label{lem:apartrinome} Let $n$ be an integer and $\gamma \in]0,1/4]$ (later defined in Lemma \ref{lem:trinome}).
Then if $$|q_h^{n+2}-q_h^{n+1}| \leq \eta \quad \textrm{and} \quad |q_h^{n+1} - q_h^{n} +h^2f_h^{n+1}|\leq \gamma \eta $$
then
$$|q_h^{n+2}-q_h^{n+1}| \leq | q_h^{n+1} - q_h^{n} +h^2f_h^{n+1} | + J |q_h^{n+1} - q_h^{n} +h^2f_h^{n+1}|^2.$$
\end{lem}

\dem By definition, for each integer $n$,
$$ q_h^{n+2} := \PPP_{C} \left[2q_h^{n+1} - q_h^{n} +h^2f_h^{n+1}\right].$$
It follows that with $z:=2q_h^{n+1} - q_h^{n} +h^2f_h^{n+1}$, $z-q_h^{n+2}$ is a proximal normal vector of $C$ at the point $q_h^{n+2}$. The hypomonotonicity property of the proximal normal cone yields
$$ \langle q_h^{n+2}-z, q_h^{n+2}-q_h^{n+1}\rangle \leq \frac{1}{2\eta} |q_h^{n+2}-q_h^{n+1}|^2 \, |z-q_h^{n+2}|.$$
It also comes
$$ |q_h^{n+2}-q_h^{n+1}|^2 - \langle z-q_h^{n+1}, q_h^{n+2}-q_h^{n+1}\rangle \leq \frac{1}{2\eta} |q_h^{n+2}-q_h^{n+1}|^2 \left(|z-q_h^{n+1}| + |q_h^{n+2}-q_h^{n+1}|\right). $$ 
Using Cauchy-Schwarz inequality, we deduce
$$  |q_h^{n+2}-q_h^{n+1}| - |z-q_h^{n+1}| \leq \frac{1}{2\eta} |q_h^{n+2}-q_h^{n+1}| \left(|z-q_h^{n+1}| + |q_h^{n+2}-q_h^{n+1}|\right). $$ 
Thanks to Lemma \ref{lem:trinome} (later proved) with $ a:=|q_h^{n+2}-q_h^{n+1}|\leq \eta$ and $b:=|z-q_h^{n+1}|\leq \gamma\eta$, it comes
$$|q_h^{n+2}-q_h^{n+1}| \leq |z-q_h^{n+1}| + J |z-q_h^{n+1}|^2.$$
\findem

\begin{lem} \label{lem:trinome}
Let $a$ and $b$ two nonnegative reals satisfying
\be{eq:trinome}  a^2+(b-2\eta) a+2\eta b \geq  0. \ee
Then for some numerical constant $\gamma\in]0,1/4[$ and $J>0$ (independent on $a,b$), we have
\be{eq:impl} \left. \begin{array}{l} 
   a \leq \eta  \\ 
 b\leq \gamma \eta 
  \end{array} \right\} \Longrightarrow a \leq b+Jb^2.
  \ee
 \end{lem}

\dem  We remark that (\ref{eq:trinome}) is a second degree polynomial function with respect to $a$, whose discriminant is given by
$$ \Delta:= (b-2\eta)^2 - 8\eta b.$$
The real $\Delta$ is nonnegative as soon as $b\leq \gamma \eta \leq \eta/4$. Then we know that (\ref{eq:trinome}) implies
$$ a\notin ]x_-,x_+[$$
with 
$$ x_\pm := \frac{-b+2\eta}{2} \pm \frac{\sqrt{\Delta}}{2}.$$
It comes
$$ \Delta = 4\eta^2 \left[1-3\frac{b}{\eta} + \frac{b^2}{4\eta^2} \right].$$
Hence,
$$ \left| \sqrt{\Delta} - 2\eta \left(1-\frac{3}{2\eta} b\right) \right| \leq M b^2$$
for some numerical constant $M=M(\gamma,\eta)$.
$$ \left|x_\pm  - \left(-\frac{b}{2}+\eta \pm\left(\eta - \frac{3b}{2} \right) \right) \right| \leq \frac{M}{2} b^2.$$
So $x_- \leq b+Mb^2/2$ and $x_+\geq 2\eta-2b-Mb^2/2$.
Furthermore, there exists a small enough constant $\gamma>0$ such that $2\gamma \eta + M( \gamma \eta)^2/2\leq \eta$, which implies with $b\leq \gamma \eta$ that $x_+\geq \eta$. Consequently, since $a\notin ]x_-,x_+[$ and $a\leq \eta$, we deduce that necessarily 
$$ a \leq x_- \leq b+\frac{M}{2} b^2,$$
which concludes the proof of (\ref{eq:impl}). \findem

\begin{lem} \label{lem:iteration} Let $J'$ be a fixed positive real. We denote by $\phi_{J'}^n$ the $n$-th iterated of 
$$\phi_{J'}:=x\mapsto x+J' x^2.$$
For every $x\geq 0$ with $J'n x \neq 1$, we have
\be{eq:lemiter} \phi_{J'}^n(x) \leq x \left( \frac{(J'nx)^{2^n-1}-1 }{J'nx-1} \right)=x \left( 1+ (J'nx) + \cdots + (J'nx)^{2^n-1}\right). \ee
\end{lem}

\dem It is obvious that $\phi_{J'}^n$ is a polynomial function of degree $2^n$ with nonnegative coefficients. 
So we know that $\phi_{J'}^{n}$ can be written as
$$ \phi_{J'}^{n}(x) =  \sum_ {k=1}^{2^{n}} a_k^{(n)} x^k.$$
We want to prove that for all $n\geq 1$ and $k\in\{1,...,2^{n}\}$
\be{eq:coeff} a_{k}^{(n)} \leq (J'n)^{k-1}.\ee
It is obvious that (\ref{eq:coeff}) holds for $n=1$ (in fact, there is equality).
Let us assume that (\ref{eq:coeff}) holds for an integer $n$ and prove it for $n+1$. Since
$$ \phi_{J'}^{n+1}(x) =  \phi_{J'}^n(x) + J' \phi_{J'}^n(x) ^2,$$
for all $k\in\{1,...,2^{n+1}\}$
\begin{align*}
 a_k^{(n+1)} & = a_k^{(n)} + J' \sum_{j=1}^{k-1} a_j^{(n)} a_{k-j}^{(n)} \\ 
 & \leq (J'n)^{k-1} + J' \sum_{j=1}^{k-1} (J'n)^{j-1} (J'n)^{k-1-j} \\
 & \leq (J'n)^{k-1} + J' (k-1) (J'n)^{k-2} \\
 & \leq (J'(n+1))^{k-1}.
\end{align*}
By agreement, $a_{k}^{(n)}$ is set equal to $0$ if $k> 2^{n}$.\\
That ends the recursive proof.
\findem

\mb We now come back to the proof of Theorem \ref{thm:appendix}.

\noindent {\bf Proof of Theorem \ref{thm:appendix}: } We define a sequence $(x^n)_n$ as follows:
$$ x^0:=|q_h^1-q_h^0|+h^2\|F\|_\infty + h\sqrt{\|F\|_\infty} \quad \textrm{and} \quad x^{n+1}=\phi_{J+1}(x^n),$$
where $J$ is introduced in Lemma \ref{lem:trinome} and functions $\phi_J$ in Lemma \ref{lem:iteration}. 
We choose $h<\min\{T_0,1\}$ small enough in order that $x_0\leq \frac{1}{2}\gamma \eta$ (this is possible since
$|q_h^1-q_h^0|\leq h|u_h^1|\leq 2h|u_0| + 2h^2\|F\|_\infty$, see Lemma \ref{lem:velocity}). \\
Now we set
$$ {\mathcal P}:=\left\{ n\leq T/h, \ |q_h^{n+1}-q_h^n|\leq \eta \ \textrm{ and } \ |q_h^n-q_h^{n-1}|+h^2|f_h^n|\leq \gamma \eta\right\}$$
and 
$$ P:= \min\left\{ n\geq 0,\  n\notin {\mathcal P}\right\}-1.$$
We want to prove that $P=N$ where $N:=T_0/h$.
Let fix $n\in \{1,...,P\}$, since $|q_h^{n+1}-q_h^n|\leq \eta$ and $|q_h^n-q_h^{n-1}|+h^2|f_h^n|\leq \gamma \eta$, thanks to Lemma \ref{lem:apartrinome} it comes
\begin{align*}
 |q_h^{n+1}-q_h^n|& \leq | q_h^{n} - q_h^{n-1} +h^2f_h^{n} | + J |q_h^{n} - q_h^{n-1} +h^2f_h^{n}|^2 \\
 & \leq \phi_{J}(| q_h^{n} - q_h^{n-1} +h^2f_h^{n} |).
\end{align*}
Consequently
\begin{align*}
 |q_h^{n+1}-q_h^n| +h^2|f_h^{n+1}| & \leq \phi_J (|q_h^{n}-q_h^{n-1}|+h^2|f_h^n|)+ h^2 \|F\|_\infty \\
 & \leq \phi_J (|q_h^{n}-q_h^{n-1}|+h^2|f_h^n|)+ (x^{0})^2 \\
 & \leq \phi_J (|q_h^{n}-q_h^{n-1}|+h^2|f_h^n|)+ (x^{n-1})^2,
 \end{align*}
where we have used that the sequence $(x^n)_n$ is increasing. Consequently, we easily deduce by iteration that for all $n=1,...,P$,
\be{eq:recur} |q_h^{n+1}-q_h^n| +h^2|f_h^{n+1}| \leq x^n.\ee
Indeed, (\ref{eq:recur}) is satisfied for $n=0$. Moreover assuming (\ref{eq:recur}) for $n-1$, we have
 $$ |q_h^{n+1}-q_h^n| +h^2|f_h^{n+1}|  \leq \phi_J (|q_h^{n}-q_h^{n-1}|+h^2|f_h^n|)+ (x^{n-1})^2 \leq \phi_{J+1} (x^{n-1}) = x^{n},$$
because $\phi_{J}$ is nondecreasing on $[0,\infty[$. \\
Thanks to Lemma \ref{lem:iteration} (with $J'=J+1$), we have for every $n\leq P$
$$ |q_h^{n+1}-q_h^n| +h^2|f_h^{n+1}|  \leq x^n \leq \left(\frac{1}{1-(J+1)nx^0}\right) x^0.$$
Moreover for all $n\leq N=T_0/h$
\begin{align*}
(J+1)n x^0 & \leq (J+1)\frac{T_0}{h} \left(|q_h^1-q_h^0|+h^2\|F\|_\infty + h\sqrt{\|F\|_\infty}\right) \\
 & \leq (J+1)T_0 \left(|u^1_h|+h\|F\|_\infty + \sqrt{\|F\|_\infty}\right) \\
 &  \leq (J+1)T_0 \left(2|u_0|+3h\|F\|_\infty + \sqrt{\|F\|_\infty}\right) \\
 & \leq \frac{1}{2},
\end{align*}
where we have used Lemma \ref{lem:velocity} to estimate $|u_h^1|$ by $2|u_0|+2h\|F\|_\infty$ and the definition of $T_0$. Thus for $n=P$, we get
$$ |q_h^{P+1}-q_h^P| +h^2|f_h^{P+1}| \leq 2 x^0\leq \gamma \eta.$$
Due to Lemma \ref{lem:velocity}, it comes
$$ |q_h^{P+2}-q_h^{P+1}| \leq 2 |q_h^{P+1}-q_h^P| + 2h^2|f_h^{P+1}| \leq 2\gamma\eta \leq \eta.$$
Finally, $P+1 \in {\mathcal P}$ and so $P=N=T_0/h$. We also conclude that
\begin{align*} 
\sup_{0\leq t_h^n \leq T_0} |u_h^n| & \leq \sup_{n\leq N} \frac{|q_h^{n}-q_h^{n-1}|}{|h|} \\
 & \leq 2 \frac{|x_0|}{h} \\
 & \leq 2\left(|u_h^1| + h \|F\|_\infty + \sqrt{\|F\|_\infty}\right) \\
 & \leq 2\left(2|u_0| + 3h \|F\|_\infty + \sqrt{\|F\|_\infty}\right),
\end{align*} 
which is uniformly bounded with respect to $h$ when $h$ goes to $0$.
\findem

\end{document}